\documentclass{gtmon_a}
\pdfoutput=1

\usepackage[all]{xy}
\usepackage{pinlabel}
\usepackage{stmaryrd}

\proceedingstitle{Proceedings of the Nishida Fest (Kinosaki 2003)}
\conferencestart{28 July 2003}
\conferenceend{8 August 2003}
\conferencename{International Conference in Homotopy Theory}
\conferencelocation{Kinosaki, Japan}

\editor{Matthew Ando}
\givenname{Matthew}
\surname{Ando}

\editor{Norihiko Minami}
\givenname{Norihiko}
\surname{Minami}

\editor{Jack Morava}
\givenname{Jack}
\surname{Morava}

\editor{W Stephen Wilson}
\givenname{W Stephen}
\surname{Wilson}

\title{Some root invariants at the prime 2}

\author{Mark Behrens}
\givenname{Mark}
\surname{Behrens}
\address{Department of Mathematics\\
Massachusetts Institute of Technology\\\newline
Cambridge MA 02139\\
USA}
\email{mbehrens@math.mit.edu}
\urladdr{}

\volumenumber{10}
\issuenumber{}
\publicationyear{2007}
\papernumber{1}
\startpage{1}
\endpage{40}

\doi{}
\MR{}
\Zbl{}

\arxivreference{math/0507183}


\received{2 August 2004}
\revised{15 May 2005}
\accepted{}
\published{29 January 2007}
\publishedonline{29 January 2007}
\proposed{}
\seconded{}
\corresponding{}
\version{}

\let\xysavmatrix\xymatrix
\def\xymatrix{\disablesubscriptcorrection\xysavmatrix}
\AtBeginDocument{\let\bar\wbar\let\tilde\wtilde\let\hat\what}

\makeop{alg}

\newcommand{\abs}[1]{\lvert #1 \rvert}

\newcommand{\br}[1]{\overline{#1}}

\newcommand{\td}[1]{\widetilde{#1}}

\newcommand{\bra}[1]{\langle #1 \rangle}
\newcommand{\ZZ}{\mathbb{Z}}
\newcommand{\QQ}{\mathbb{Q}}
\newcommand{\RR}{\mathbb{R}}
\newcommand{\CC}{\mathbb{C}}
\newcommand{\FF}{\mathbb{F}}
\newcommand{\capeq}{\overset{\scriptscriptstyle{\cap}}{=}}
\newcommand{\dotin}{\overset{\scriptstyle{\centerdot}}{\in}}

\makeatletter
\def\cnewtheorem#1[#2]#3{\newtheorem{#1}{#3}[section]
\expandafter\let\csname c@#1\endcsname\c@thm}

\newtheorem{thm}{Theorem}[section]
\cnewtheorem{cor}[thm]{Corollary}
\cnewtheorem{lem}[thm]{Lemma}
\cnewtheorem{prop}[thm]{Proposition}
\theoremstyle{remark}
\cnewtheorem{rmk}[thm]{Remark}
\cnewtheorem{defn}[thm]{Definition}

\newtheorem*{conventions}{Conventions}
\cnewtheorem{procedure}[thm]{Procedure}

\makeatother  

\DeclareMathOperator{\ext}{Ext}

\newcommand{\lk}{\mathit{lk}}

\begin{document}

\begin{htmlabstract}
The first part of this paper consists of lecture notes which summarize the
machinery of filtered root invariants.  A conceptual notion of "homotopy
Greek letter element" is also introduced, and evidence is presented that
it may be related to the root invariant.  In the second part we compute
some low dimensional root invariants of v<sub>1</sub>&ndash;periodic
elements at the prime 2.
\end{htmlabstract}

\begin{abstract}
The first part of this paper consists of lecture notes which summarize the
machinery of filtered root invariants.  A conceptual notion of ``homotopy
Greek letter element'' is also introduced, and evidence is presented that
it may be related to the root invariant.  In the second part we compute some
low dimensional root invariants of $v_1$--periodic elements at the prime
$2$.
\end{abstract}

\maketitle


\section{Introduction}

This paper consists of two parts.  The first part consists of the lecture
notes of a series of talks on the root invariant given by the author at a 
workshop held at the Nagoya Institute of Technology.  
The second part is a detailed computation, using the
methods of \fullref{part:lec}, of some low dimensional root invariants at
the prime $2$.  More detailed descriptions of the contents of these parts
are given at the beginning of each part.

\begin{tabular}{l}
Part I\qua Lectures on root invariants \\
\quad 2\qua The chromatic filtration \\
\quad 3\qua Greek letter elements \\
\quad 4\qua The root invariant \\
\quad 5\qua Filtered root invariants \\
\quad 6\qua Some theorems \\
Part II\qua $2$--primary calculations \\
\quad 7\qua The indeterminacy spectral sequence \\
\quad 8\qua $BP$ filtered root invariants of $2^k$ \\
\quad 9\qua The first two $BP$--filtered root invariants of $\alpha _{i/j}$ \\
\quad 10\qua Higher $BP$ and $H\mathbb{F}_2$--filtered root invariants of some 
$v_1$--periodic elements \\
\quad 11\qua Homotopy root invariants of some $v_1$--periodic elements
\end{tabular}

The author would like to express his appreciation of the contributions of
Goro Nishida to the field of homotopy theory.  The author would also like
to thank the organizers of the conference for the unique opportunity to
meet Japanese mathematicians.  Jack Morava encouraged the author 
to submit these
lecture notes as a means of communicating the ideas of \cite{Behrensroot}
without all of the technical details.  Haynes Miller provided some comments
on the homotopy Greek letter construction, and Mike Hopkins clued the
author into the existence of Mahowald's useful 
paper \cite{MahowaldDescriptive}.  
W-H~Lin explained to the author how to show that a certain element of
the Adams spectral sequence for $P_{-23}$ was a permanent cycle.
The author is grateful to the referee for discovering a mistake in a
previous version of \fullref{corA}, and to R\,R~Bruner for pointing
out some typographical errors.
The computations of
\fullref{part:calc} began as part of the author's thesis, which was
completed under the guidance of Peter May at the University
of Chicago.  The
author is, here and elsewhere, heavily influenced by the work of (and
discussions with) Mark Mahowald.  Finally, the author would like to extend
his heartfelt gratitude to Norihiko Minami, both for organizing a very engaging
workshop at the Nagoya Institute of Technology, and for his role as a
mentor, introducing the author to
the field of homotopy theory as an undergraduate at the University of
Alabama.  
The author is supported by the NSF.

\part{Lectures on root invariants}
\setobjecttype{Part}
\label{part:lec}

Throughout this paper
we are always working in the $p$--local stable homotopy category for $p$
a fixed prime number.  In this first section we will summarize the chromatic
filtration and the machinery of filtered root invariants.  A very detailed
treatment of this theory appeared in \cite{Behrensroot}.  Our intention
here is to ignore many of the subtleties, sometimes to the point of
omitting or simplifying hypotheses and ignoring indeterminacy, to
communicate to the reader the underlying ideas.  Our hope is that the
reader will be able to use this overview as a motivation to drudge through
the more precise treatment of \cite{Behrensroot}.

In \fullref{sec:chromatic} we review the chromatic filtration of the
stable stems.  In \fullref{sec:Greek} we review the Greek letter
construction of Miller, Ravenel, and Wilson.  The Greek letter elements are
distinguished chromatic families of elements in the $E_2$--term 
$\ext_{BP_*BP}(BP_*,BP_*)$ of the Adams--Novikov spectral sequence (ANSS).
These elements need not be non-trivial permanent cycles in the
ANSS.  We introduce the notion of a
homotopy Greek letter element to remedy this.  In \fullref{sec:rootdef}
we define the root invariant and recall some computational examples that
occur throughout the literature.  The interesting thing is that, at least
for the limited number of root invariants we know, it seems to be the case
that the root invariant has a tendency 
to take $v_n$--periodic homotopy Greek letter elements to $v_{n+1}$--periodic 
homotopy Greek letter elements.  In \fullref{sec:filtrootdef} we define
filtered root invariants.  In \fullref{sec:results} we summarize the
main theorems that make the filtered root invariants compute root
invariants.

\begin{conventions}
We shall be using the following abbreviations.
\begin{description}
\item[ASS] The classical Adams spectral sequence based on $H\FF_p$.
\item[ANSS] The Adams--Novikov spectral sequence based on $BP$.
\item[$E$--ASS] The generalized Adams spectral sequence based on $E$.
\item[AHSS] The Atiyah--Hirzebruch spectral sequence.
\item[$\ext(E_*X)$] For nice spectra $E$, the comodule $\ext$
  group $\ext_{E_*E}(E_*, E_*X)$.
\item[AAHSS] The algebraic Atiyah--Hirzebruch spectral sequence that
  computes the group $\ext(E_*X)$ using the cellular filtration on $X$.
\item[Modified AAHSS] The AAHSS for $\ext(BP_*P^\infty)$ arising from the
  filtration of $P^\infty$ by Moore spectra.
\end{description}
Often we shall place a bar over the name of a permanent cycle in an Adams
spectral sequence to denote an element of homotopy that it detects.  We
shall place dots over binary relations to indicate that they only hold up
to multiplication by a unit in $\ZZ_{(p)}$.  For instance, we shall write
$x \doteq y$ if $x = \alpha \cdot y$ for some $\alpha \in
\ZZ_{(p)}^\times$.  We shall use $\smash{\capeq}$ to indicate that two quantities
are equal modulo some indeterminacy group.  
We shall always use the Hazewinkel generators of
$BP_*$.
\end{conventions}

\section{The chromatic filtration}\label{sec:chromatic}

We shall first describe the chromatic filtration on the stable homotopy
groups of spheres.  What we are describing is referred to as the
``geometric chromatic filtration'' by Ravenel in \cite{Ravenelorange}.
We first need to discuss type--$n$ complexes and $v_n$--self maps.

Let $K(n)$ be $n^\mathrm{th}$ Morava $K$--theory, with coefficient ring
$$ K(n)_* = \FF_p[v_n, v_n^{-1}]. $$
Here $v_n$ has degree $2(p^n-1)$.  By $K(0)$ we shall mean the rational
Eilenberg--MacLane spectrum $H\QQ$, and by $v_0$ we shall mean $p$.

If $X$ is a finite complex, it is said to be \emph{type--$n$} if $K(n-1)_*X
= 0$ and $K(n)_*X \ne 0$.  It is a consequence of the Landweber filtration
theorem (see Landweber \cite{Landweber} and Ravenel \cite{Ravenelorange})
that the condition $K(n-1)_*X = 0$ implies that $K(m)_*X = 0$ for $m
\le n-1$.

A self map $v\co  \Sigma^{N}X \rightarrow X$ is said to be a $v_n$--self map if
it induces an isomorphism on $K(n)$--homology
$$ v_*\co  K(n)_*\Sigma^NX \rightarrow K(n)_*X. $$
If $v_*$ induces,
up to an element in
$\smash{\FF_p^\times}$, multiplication by $v_n^k$ for some $k$, we shall say that
$X$
\emph{has $v_n^k$ multiplication}. 
By \cite[Lemma~6.1.1]{Ravenelorange}, if $v$ is a $v_n$--self map, 
there is some power of $v$ which induces $v_n^k$ multiplication.   
If $v$ induces $v_n^k$--multiplication, we shall often denote the map
$v_n^k$.
This practice is slightly objectionable
because a complex can have many different $v$ so that $v_*$ is $v_n^k$, but
there is some consolation in that 
the Devinatz--Hopkins--Smith Nilpotence Theorem may be used 
to show that any two such maps are equal
after a finite number of iterates.

The following important theorem was conjectured by Ravenel
\cite{Ravenel84} and
proved by Hopkins and Smith \cite{NilpI} using the nilpotence theorem
\cite{NilpII}.

\begin{thm}[Hopkins--Smith Periodicity Theorem]
If $X$ is type--$n$, then $X$ possesses a $v_n$--self map.
\end{thm}

We will now define the chromatic filtration of an element 
$\alpha \in \pi_n(S)$.  We shall refer
to the following diagram.
$$
\xymatrix{ S^n \ar[r]_{p^{k_0}} & 
S^n \ar[d] \ar[r]^\alpha & 
S^0
\\
\Sigma^{N_1} M(p^{k_0}) \ar[r]_{v_1^{k_1}} &
M(p^{k_0}) \ar@{.>}[ru]^{\alpha_1} \ar[d]
\\
\Sigma^{N_2} M(p^{k_0}, v_1^{k_1}) \ar[r]_-{v_2^{k_2}} &
M(p^{k_0}, v_1^{k_1}) \ar@{.>}[ruu]_{\alpha_2} \ar[d] 
\\
& \vdots
}
$$

\begin{description}
\item[$v_0$--periodic]
If $\alpha\circ p^k$ is non-zero for every $k$, then $\alpha$ is said to be
\emph{$v_0$--periodic}.

\item[$v_1$--periodic]
Otherwise, $\alpha$ is $v_0$--torsion, and there is some
$k_0$ such that $\alpha \circ p^{k_0} = 0$.  Let $M(p^{k_0})$ denote the
cofiber of $p^{k_0}$.  Then there exists a lift $\alpha_1$ of $\alpha$ to
$M(p^{k_0})$.  The complex $M(p^{k_0})$ is type--$1$, and thus has $v_1^m$
multiplication for some $m$.  If $\alpha_1\circ \smash{v_1^{mk}}$ is non-zero for
every $k$, then $\alpha$ is said to be \emph{$v_1$--periodic}.

\item[$v_2$--periodic]
Otherwise, $\alpha$ is $v_1$--torsion, and there is some
$k_1$ with $\alpha_1{\circ}v_1^{k_1}{=}0$.  Let $M(\smash{p^{k_0}},
\smash{v_1^{k_1}})$ denote the cofiber of $\smash{v_1^{k_1}}$.  Then there
exists a lift $\alpha_2$ of $\alpha_1$ to
$M(\smash{p^{k_0}}, \smash{v_1^{k_1}})$.  The complex
$M(\smash{p^{k_0}},\smash{v_1^{k_1}})$ is type--$1$, and thus has $v_2^m$
multiplication for some $m$.  If $\alpha_2\circ \smash{v_2^{mk}}$ is non-zero for
every $k$, then $\alpha$ is said to be \emph{$v_2$--periodic}.

\item[$v_3$--periodic]
Otherwise, $\alpha$ is said to be $v_1$--torsion, and there is some $k_2$ so 
that $\alpha_2 \circ \smash{v_2^{k_2}} = 0$.  The process continues.

\end{description}

In this way, we have defined a decreasing filtration
$$ \pi_*(S) \supseteq \{v_0-\mathit{torsion}\} \supseteq
\{v_1-\mathit{torsion}\}
\supseteq \{v_2-\mathit{torsion}\} \supseteq \cdots
$$
which is the chromatic filtration.  

It is not clear that the chromatic filtration
is independent of the sequence of lifts.
The (geometric) chromatic filtration may be more succinctly described by
means of finite localization (see Mahowald--Sadofsky
\cite{MahowaldSadofsky}), and from this
perspective it is clear that the chromatic filtration is well defined.  
The \emph{finite localization} functor 
$$ L^f_{E(n)}\co  \mathit{Spectra} \rightarrow \mathit{Spectra} $$
is initial amongst endofunctors that kill finite $E(n)$--acyclic spectra.
The finitely localized sphere $L_{E(n)}^fS$ may be explicitly described as a
colimit of finite spectra, and in this manner one finds that the
$v_{n}$--torsion elements of $\pi_*(S)$ are precisely those that make up the
kernel of the map
$$ \pi_*(S) \rightarrow \pi_*(L_{E(n)}^fS). $$
There are fiber sequences
$$ M_n^fS = v_n^{-1}M(p^\infty, \ldots , v_{n-1}^\infty) \rightarrow L_{n}^fS
\rightarrow L_{n-1}^fS. $$

\begin{rmk}
In \cite{Ravenelorange}, Ravenel discusses a different filtration which he
calls the ``algebraic chromatic filtration,'' which is what is more
commonly meant by the chromatic filtration these days.  The $n^\mathrm{th}$
filtration is the kernel of the localization map
$$ \pi_*(S) \rightarrow \pi_*(L_{E(n)}S) $$
where $E(n)$ is the Johnson--Wilson spectrum with coefficient ring
$$ E(n)_* = \ZZ_{(p)}[v_1, v_2, \ldots, v_n, v_n^{-1} ].$$
If the telescope conjecture is true, than the ``geometric'' and
``algebraic'' chromatic filtrations agree.  However, it has been the case
in the past decade that the more people have thought about the telescope
conjecture, the more they have believed it to be false
(see Mahowald--Ravenel--Shick \cite{MahowaldRavenelShick}).
\end{rmk}

\section{Greek letter elements}\label{sec:Greek}

We shall now outline a standard method of constructing $v_n$--periodic
elements of the stable stems called the \emph{Greek letter construction}.
Suppose that the generalized Moore spectrum $M(I)$ exists, where $I$ is the
ideal $(\smash{p^{i_0}, v_1^{i_1}, \ldots, v_{n-1}^{i_{n-1}}}) \subset BP_*$, 
and assume that $M(I)$ has $v_n^k$--multiplication.  
The spectrum $M(I)$ is a finite complex of dimension
$$ d = n+i_1\abs{v_1} + \cdots + i_{n-1}\abs{v_{n-1}}.$$
Then we can form the composite
$$ \alpha^{(n)}_{\lk/i_{n-1},\ldots,i_{0}} \co 
S^{\lk\abs{v_n}-d} \xrightarrow{\iota} \Sigma^{\lk\abs{v_n}-d}M(I)
\xrightarrow{v_n^{\lk}} \Sigma^{-d}M(I) \xrightarrow{\nu} S^0 $$
where $\iota$ is the inclusion of the bottom cell, $\nu$ is the projection
onto the top cell, and 
$\alpha^{(n)}$ is the $n^\mathrm{th}$ letter in the Greek alphabet
$\alpha, \beta, \gamma, \delta, \ldots $.
Miller, Ravenel, and Wilson in \cite{MillerRavenelWilson} described Greek
letter elements in $\ext(BP_*)$, and the Greek letter elements of
$\pi_*(S)$ that we have described are detected by their elements in the
Adams--Novikov spectral sequence (ANSS).  We shall refer to the Greek letter
elements in $\ext(BP_*)$ as ``algebraic Greek letter elements.''

We give a different interpretation of the Greek letter construction with an
eye towards generalization.  The existence of $v_n^k$ multiplication on
$M(I)$ gives homotopy elements
$$ v_n^{\lk} \in \pi_{\lk\abs{v_n}}(M(I)) $$
detected by $v_n^{\lk}$ in $BP$--homology.  Fix a (minimal) 
cellular decomposition of
$M(I)$.  Consider the
Atiyah--Hirzebruch spectral sequence (AHSS)
$$ E_1^{n,i} = \bigoplus_{\text{$n$--cells in $M(I)$}} \pi_i(S^n)
\Rightarrow \pi_i(M(I)). $$
Suppose that the element $\smash{\alpha^{(n)}_{\lk/i_{n-1}, \ldots, i_0}}$ is
non-trivial.  Then in the AHSS, $v_n^{\lk}\in \pi_*(M(I))$ is detected on
the top cell by $\smash{\alpha^{(n)}_{\lk/i_{n-1}, \ldots, i_{0}}}$.

But how does one define $\smash{\alpha^{(n)}_{k/i_{n-1},\ldots, i_{0}}}$
if the appropriate $M(I)$ does not exist?  Or if $M(I)$ does not
have $v_n^k$ multiplication?  What do we do if the homotopy element
$\smash{\alpha^{(n)}_{k/i_{n-1},\ldots, i_{0}}}$ turns out to be trivial?
We give a ``definition'' of \emph{homotopy Greek letter elements} to be
the homotopy replacement of the Greek letter element when any of the
above calamities befalls us.  The author does not believe this is the
right way to define these elements, but has no better ideas.

\begin{defn}[Homotopy Greek letter elements]
Suppose $X$ is a type--$n$ $p$--local finite complex for which $BP_*X$ is
free module over $BP_*/I$, for 
$I = (p^{i_0}\!, v_1^{i_1}\!, \ldots, v_{n-1}^{i_{n-1}})$.
Suppose that $X$ has $v_n^k$--multiplication.  Then we define the homotopy
Greek letter element $\smash{(\alpha^{(n)})^h_{k/i_{n-1}, \ldots, i_0}}$
to be the element of $\pi_*(S)$ which detects $v_n^k \in \pi_*(X)$
in the $E_1$--term of the AHSS.
\end{defn}

\begin{rmk}
This definition is full of flaws.  Different
choices of $X$, or even different choices of detecting element in the AHSS,
could yield some ambiguity in the definition of 
$\smash{(\alpha^{(n)})^h_{k/i_{n-1}, \ldots, i_0}}$.  
There is no reason to believe that these homotopy Greek letter elements are
of chromatic filtration $n$.
In fact, the stem in which
the element appears is even ambiguous.  We do point out that there is
already some ambiguity in the standard definition of the Greek letter
elements --- there can exist many complexes with the same $BP$--homology as
$M(I)$, and the choice of $v_n$--self map is not unique.  The $v_n$--self map 
is, however, unique after a finite number of iterations (see
Hopkins--Smith \cite{NilpII}).  
\end{rmk}

\begin{rmk}
Mahowald and Ravenel \cite{MahowaldRavenel}
propose a different definition for ``homotopy Greek
letter elements'' using iterated root invariants.  Their definition
suggests that we should be defining $(\alpha^{(n)})^h_i$ as
$$ (\alpha^{(n)})^h_i = R^n(p^i) $$
This notion suffers the same sorts of indeterminacy issues that our notion of
homotopy Greek letter element suffers.  
\end{rmk}

We give some examples to illustrate the sorts of phenomena that the reader
should expect at bad primes.

Let $p = 2$, and consider chromatic level $n = 1$.  The complex $M(2)$ only
has $v_1^4$ multiplication (see Adams \cite{Adams}), giving us the Greek letter
elements $\alpha_{4k} \in \pi_{8k-1}(S)$ (these are elements of order $2$ in
the image of $J$).  However the complex $X = M(2) \wedge C(\eta)$, where
$C(\eta)$ is the cofiber of $\eta$, has $v_1$--multiplication (see Mahowald
\cite{Mahowald}).  
Using this
complex, we get the following homotopy Greek letter elements.  These are
precisely the elements on the edge of the ASS vanishing line, and, we
believe, quite worthy of the designation ``Greek letter element.''

\begin{center}
\begin{tabular}{|c|c|}
\hline
$n \pmod 4$ & $\alpha^h_{n}$ (ANSS name) \\
\hline
$0$ & $\alpha_{n}$ \\
$1$ & $\alpha_{n}$ \\
$2$ & $\alpha_{n-1}\alpha_1$ \\
$3$ & $\alpha_{n-2}\alpha_1^2$ \\
\hline
\end{tabular}
\end{center}

Likewise, we list some low dimensional homotopy Greek letters for $p=3$ and
chromatic level $n = 2$.  The complex $M(3,v_1) = V(1)$ exists, but only has
$v_2^9$ multiplication (see Behrens--Pemmaraju \cite{BehrensPemmaraju}).  Thus we are only able to
define $\beta_{9t}$ using the conventional methods.  However the complex
$$ X = V(1)\wedge(S^0 \cup_{\alpha_1} e^4 \cup_{\alpha_1} e^8) $$
can be shown to have $v_2$--multiplication, and using this complex we find
the following homotopy Greek letter elements.

\begin{center}
\begin{tabular}{|c|c|}
\hline 
Greek name & Adams--Novikov name \\
\hline
$\beta_1^h$ & $\beta_1$ \\
$\beta_2^h$ & $\beta_1^2\alpha_1$ \\
$\beta_3^h$ & $\beta_3$ \\
$\beta_4^h$ & $\beta_1^5$ \\
$\beta_5^h$ & $\beta_5$ (*) \\
$\beta_6^h$ & $\beta_6$ (*)\\
\hline
\end{tabular}
\end{center}
(*) Tentative calculation

The reader will note that although the algebraic Greek letter element
$\beta_2$ exists in $\ext(BP_*)$ and is a non-trivial permanent cycle in
the ANSS, it does not agree with the homotopy Greek letter element
$\beta_2^h$.

We shall present evidence that these homotopy Greek letter elements may 
behave nicely with respect to root invariants.

\section{The root invariant}\label{sec:rootdef}

Mahowald defined an invariant called the \emph{root invariant} that takes an
element $\alpha$ in the stable stems and outputs another element
$R(\alpha)$ in the stable stems.  
$$ R\co  \pi_*(S) \leadsto \pi_*(S) $$
Our reason for using the wavy arrow ``$\leadsto$'' is that $R$ is not a
well defined map, but has indeterminacy, much in the way that Toda brackets
do.  $R(\alpha)$ is actually a coset, but in this first part, we shall
often ignore this indeterminacy to clarify the exposition.
In the literature this invariant is
sometimes called the ``Mahowald invariant,'' with good reason.

In this section we shall define the root invariant.
We shall then summarize some of the computations of root invariants that
appear in the literature.

We first need to define stunted projective spectra.  We first assume that
we are working at the prime $p = 2$.  Let $\xi$ be the canonical line
bundle over $\RR P^n$.  Then the Thom space may be identified
(see Bruner--May--McClure--Steinberger \cite[V.2.14]{Hinfty}) as
$$ (\RR P^n)^{s\xi} \cong \RR P^{n+s}/\RR P^{s-1}. $$
We may allow $s$ to be negative in the above definition if we use Thom
\emph{spectra} instead of Thom spaces.  This motivates the definition of
the spectrum
$$ P^{n+s}_s = (\RR P^n)^{s\xi} $$
for any integer $s$ and any non-negative integer $n$.  This spectrum has
one cell in each degree in the interval $[s,n+s]$.  

At an odd prime we can replace $\RR P^n = B\Sigma_2$ with the classifying
space $B\Sigma_p$.  This complex only has cells in degrees congruent to
$0, -1 \pmod {2(p-1)}$.  We shall also refer to the resulting spectra as
$P_s^{n+s}$. 

We can take the colimit over $n$ to obtain the spectrum $P^\infty_s$.
Taking the homotopy inverse limit of these spectra over $s$ yields a
spectrum $P_{-\infty}^\infty$.  The inclusion of the $-1$--cell extends to a
map
$$ S^{-1} \xrightarrow{l} P_{-\infty}^\infty. $$
For $p = 2$, 
Lin \cite{Lin} proved the following remarkable theorem. The theorem was 
conjectured by 
Mahowald, and is equivalent to the Segal conjecture for the group $\ZZ/2$.
The odd primary version was proved by Gunawardena \cite{Gunawardena}.

\begin{thm}
The map $l\co  S^{-1} \rightarrow P_{-\infty}^\infty$ is equivalent to the
$p$--completion of $S^{-1}$.
\end{thm}

This theorem makes the following definition possible.

\begin{defn}[Root invariant]
Let $\alpha$ be an element of  
$\pi_t(S^0)$.  The \emph{root invariant} 
of $\alpha$ is the coset 
of all dotted
arrows making the following diagram commute.
$$ \xymatrix{
S^{t-1} \ar@{.>}[r] \ar[d]_\alpha & S^{-N} \ar[dd]  
\\
S^{-1} \ar[d]_l
\\
P_{-\infty}^\infty \ar[r] & P_{-N}
}$$
This coset is denoted $R(\alpha)$.  Here $N$ is chosen to be minimal
such that the composite $S^{t-1} \rightarrow P_{-N}$ 
is non-trivial.
\end{defn}

One way to think of the root invariant is that it is the coset $R(\alpha)$
of elements in the $E_1$--term of the AHSS
$$ E_1^{k,n} = \pi_k(S^n) \Rightarrow \pi_k(P_{-\infty}^\infty) =
\pi_{k}(S^{-1}_2) $$
that detects $\alpha$.

The root invariant is interesting for two reasons:
\begin{enumerate}
\item Elements which are root invariants behave quite differently in the
EHP sequence as opposed to elements which are not root invariants
(see Mahowald--Ravenel \cite{MahowaldRavenel}).

\item The root invariant appears to generically take
things in chromatic filtration $n$ to things in chromatic filtration $n+1$.
\end{enumerate}

Our purpose in this paper is to concentrate on the latter.
For instance, we give the following sampler of results:

\begin{itemize}

\item For $p \ge 3$ we have $R(p^i) = \alpha_i$ \cite{MahowaldRavenel}.

\item For $p \ge 5$ we have $R(\alpha_i) = \beta_i$ and $R(\alpha_{p/2}) =
\beta_{p/2}$ (see \cite{MahowaldRavenel} and Sadofsky \cite{Sadofsky}).

\item For $p = 2$ we have $R(2^i) = \alpha_i^h$ (see
\cite{MahowaldRavenel}, Johnson \cite{Johnson}).

\item For $p = 3$ we have $R(\alpha_{i}) = \beta_{i}^h$ 
for $i \le 6$, and we have $R(\alpha_{i}) = \beta_i$ for for $i \equiv
0,1,5 \pmod 9$ (see Behrens \cite{Behrensroot}). 
\end{itemize}

A conjecture that the root invariant increases chromatic filtration 
appears in Mahowald--Ravenel \cite{MahowaldRavenelglobal}.
However, we warn the reader that the conjecture takes some time to begin
working.  For instance, at $p=2$, $R(\eta) = \nu$, and $R(\nu) =
\sigma$, and $\eta$, $\nu$, and $\sigma$ are all $v_1$--periodic elements.

\section{Filtered root invariants}\label{sec:filtrootdef}

Let $E$ be a ring-spectrum for which the $E$--Adams spectral sequence
converges.
In \cite{Behrensroot}, the author investigated a series of approximations
to the root invariant which live in the $E_1$--term of the $E$--Adams
spectral sequence called \emph{filtered root invariants}.  
$$ R^{[k]}_{E}\co  \pi_*(S) \leadsto E_1^{k,*} $$
We shall give a
brief outline of their definition, but refer the reader to
\cite{Behrensroot} for a completely detailed treatment.

Let $\wbar{E}$ be the fiber of the unit $S \rightarrow E$.
The $E$--Adams resolution of the sphere is given by
$$ \xymatrix{
S \ar@{=}[r] & W_0 \ar[d] & W_1 \ar[d] \ar[l] & W_2 \ar[d] 
  \ar[l] & W_3 \ar[d] \ar[l] & \cdots \ar[l] \\
& Y_0 & Y_1 & Y_2 & 
Y_3
} $$
where $W_k = \wbar{E}^{(k)}$ and $Y_k = E \wedge \wbar{E}^{(k)}$.  
The skeletal filtration of $P_{-\infty}^\infty$ is
given by
$$
\xymatrix{
\cdots \ar[r] & 
P_{-\infty}^N \ar[d] \ar[r] & 
P_{-\infty}^{N+1} \ar[d] \ar[r] & 
P_{-\infty}^{N+2} \ar[d] \ar[r] & 
\cdots 
\\
& S^N & S^{N+1} & S^{N+2} & 
}
$$

We wish to mix the two filtrations.  We may regard $P_{-\infty}^\infty$ as
being a bifiltered object, with $(k,N)$--bifiltration given by
$$ W_k(P^N) = (W_{k}\wedge P^{N})_{-\infty} $$
where we take the homotopy limit after smashing with $W_k$.
We may pictorially represent this bifiltration by a region of the Cartesian
plane where we let the $x$--axis represent the Adams filtration and the
$y$--axis represent the skeletal filtration.

\begin{center}
\labellist\small
\pinlabel {Adams filtration} [t] at 334 14
\pinlabel {\rotatebox{90}{\parbox{34pt}{Skeletal filtration}}} [r] at 29 157
\pinlabel {$N$} [r] at 29 122
\pinlabel {$k$} [t] at 119 14
\pinlabel {$W_k\bigl(P^n\bigr)$} [l] at 308 68
\endlabellist
\includegraphics[width=4in]{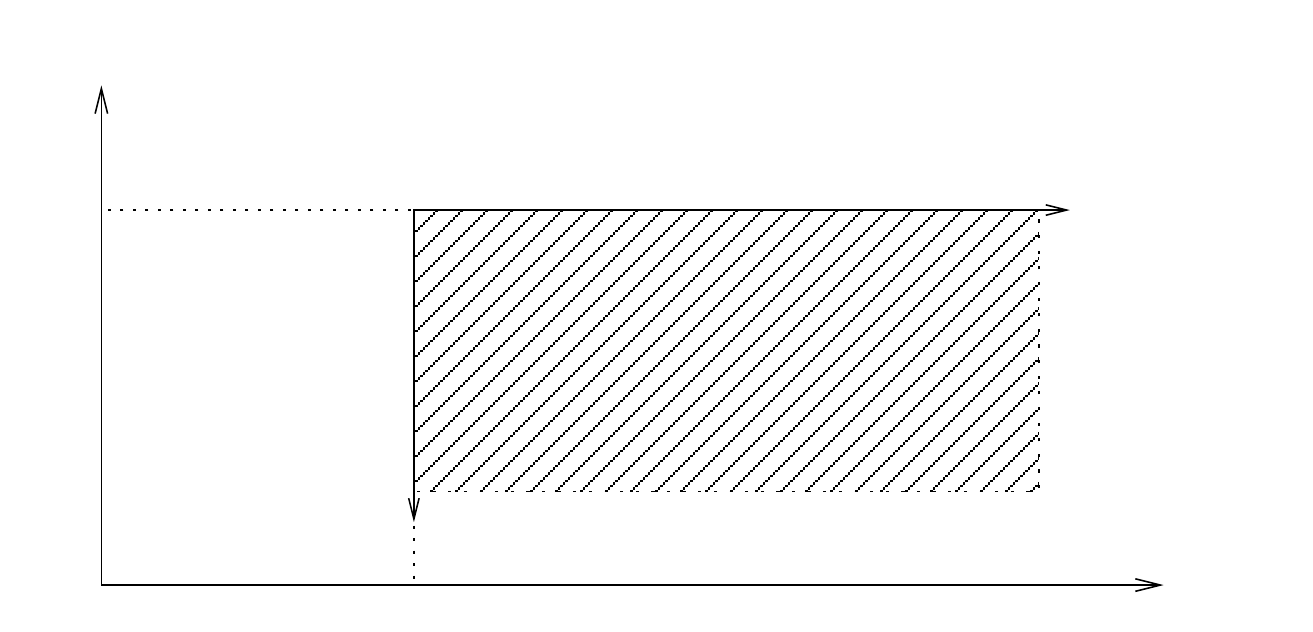}
\end{center}

The spectra $W_k$ may be replaced by weakly equivalent approximations 
so that for every $k$ the maps $W_{k+1} \rightarrow W_{k}$ 
are inclusions of
subcomplexes.  We then have that for $k_1 \ge k_2$ and $N_1 \le N_2$, the
bifiltration $W_{k_1}(P^{N_1})$ is a subcomplex of 
$W_{k_2}(P^{N_2})$.  We shall consider spectra which are
unions of these bifiltrations, which appeared in \cite{Behrensroot} as
``filtered Tate spectra.''  Given sequences 
\begin{align*}
I & = \{ k_1 < k_2 < \cdots < k_l \} \\
J & = \{ N_1 < N_2 < \cdots < N_l \}
\end{align*}
with $k_i \ge 0$,
we define the filtered Tate spectrum as the union
$$ W_I(P^J) = \bigcup_i W_{k_i}(P^{N_i}). $$
A picture of the bifiltrations that compose this spectrum is given below.

\begin{center}
\labellist\small
\pinlabel {\rotatebox{90}{\parbox{35pt}{Skeletal filtration}}} [r] at 26 217
\pinlabel {Adams filtration} [t] at 333 19
\pinlabel {$N_1$} [r] at 27 91
\pinlabel {$N_2$} [r] at 27 123
\pinlabel {$N_3$} [r] at 27 176
\pinlabel {$k_1$} [t] at 45 19
\pinlabel {$k_2$} [t] at 112 19
\pinlabel {$k_3$} [t] at 202 19
\pinlabel {\scriptsize $W_I\bigl(P^J\bigr)$}  at 260 108
\endlabellist
\includegraphics[width=4.5in]{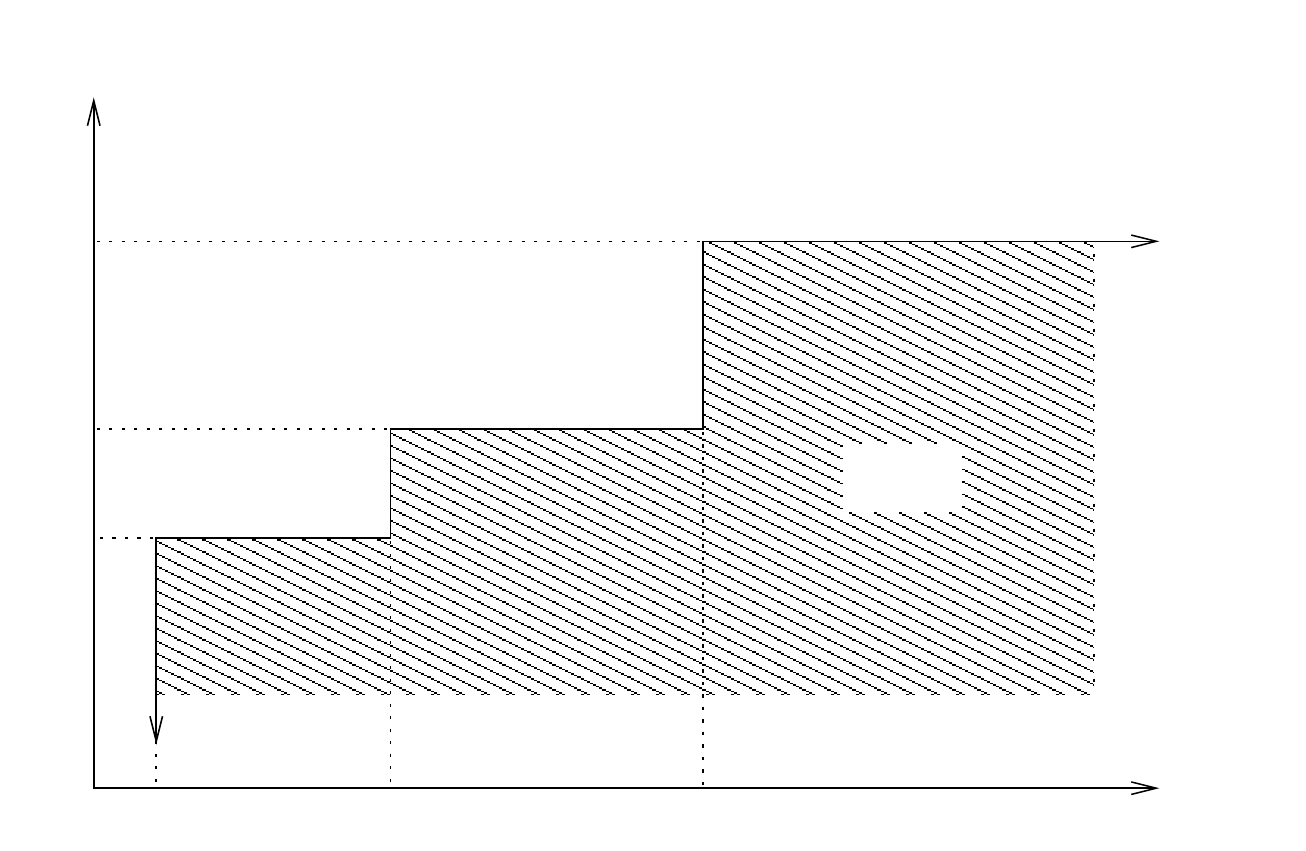}
\end{center}

For $1 \le i \le l$, there are natural projection maps
$$ p_{i}\co  W_I(P^J) \rightarrow Y_{k_i}\wedge S^{N_i}.$$
We shall now define the filtered root invariants.

\begin{defn}
Let $\alpha$ be an element of $\pi_t(S)$, with image $l(\alpha) \in
\pi_{t-1}(P_{-\infty}^\infty)$.
Choose a multi-index $(I,J)$ where $I = (k_1, k_2, \ldots)$ and $J = (N_1,
N_2, \ldots)$ so 
that the filtered Tate spectrum $W_I(P^J)$ is initial amongst the Tate
spectra $W_K(P^L)$ so that $l(\alpha)$ is in the image of the map
$$ \pi_{t-1}(W_K(P^L)) \rightarrow \pi_{t-1}(P_{-\infty}^\infty). $$
(This initial multi-index is not unique with this property, 
but in \cite{Behrensroot} we give a
convention for choosing a unique preferred initial multi-index.)
Let $\wtilde{\alpha}$ be a lift of $l(\alpha)$ to $\pi_{t-1}(W_I(P^J))$.  Then
the $k_i^\mathrm{th}$ $E$--filtered root invariant is given by
$$ R^{[k_i]}_{E}(\alpha) = 
p_i(\wtilde{\alpha}) \in \pi_{t-1}(Y_{k_i} \wedge S^{N_i}). $$
We shall refer to $(I,J)$ as the $E$--bifiltration of $\alpha$.
\end{defn}

The $k_i^\mathrm{th}$ 
filtered root invariant thus lives in the $E_1$--term of the $E$--ASS for
the sphere.  It should be regarded as an approximation to the root
invariant in $E$--Adams filtration $k_i$.  There is indeterminacy in this
invariant given by the various choices of lifts $\wtilde{\alpha}$.

\section{Some theorems}\label{sec:results}

We shall now outline the manner in which filtered root invariants may be
used to compute homotopy root invariants.  The statements of
these theorems appeared in \cite[Section~5]{Behrensroot}, with proofs
appearing in Section~6.  The theorems as stated in \cite{Behrensroot} are
rather difficult to conceptualize due to the complicated hypotheses and the
nature of the indeterminacy.  The statements we give below are imprecise,
but easier to read and understand.  
Throughout this section, let $\alpha$
be an element of $\pi_t(S)$ of $E$--bifiltration $(I,J)$, where $I =
(k_i)$ and $J = (-N_i)$.  Our first theorem tells us how to determine if
a filtered root invariant detects the homotopy root invariant in the
$E$--ASS.

\begin{thm}{\rm\cite[Theorem~5.1]{Behrensroot}}\label{thmA}\qua
Suppose that $R^{[k_i]}_E(\alpha)$ contains a permanent cycle $\beta$. Then
there exists an element
$\wbar{\beta} \in \pi_*(S)$ which $\beta$ detects in the $E$--ASS such that 
the following diagram commutes up to elements of $E$--Adams
filtration greater than or equal to $k_{i}+1$.
$$\xymatrix{
S^{t-1} \ar[d]_{\alpha} \ar[r]_-{\wbar{\beta}} & S^{-N_i} \ar[dd]  \\
S^{-1} \ar[d]_l \\
P_{-\infty}^\infty \ar[r] & P_{-N_i}
}$$
\end{thm}

We present a practical reinterpretation of this theorem.  This essentially
appears in \cite{Behrensroot} as Procedure~9.1.

\begin{cor}\label{corA}
Let $\beta$ be as in \fullref{thmA}.  Then in order for $\beta$ to
detect the homotopy root invariant in the $E$--ASS, it is sufficient to
check the following two things.
\begin{enumerate}
\item No element $\gamma \in \pi_{t-1}(P_{-N_i})$ of $E$--Adams filtration 
greater that 
$k_i$ can detect the root invariant of $\alpha$ in $P_{-N_i+1}$.

\item The image of the element $\wbar{\beta}$ under the inclusion of the
bottom cell
$$ \pi_{t-1}(S^{-N_i}) \rightarrow \pi_{t-1}(P_{-N_i})$$
is non-trivial.
\end{enumerate}
\end{cor}

For our next set of theorems we shall need to introduce a variant of the 
Toda bracket construction.  Let
$K$ be a finite CW complex with a single cell in top dimension $n$ and
bottom dimension $0$.  There is an inclusion map 
$$
\iota \co  S^0 \rightarrow K
$$
and the $n$--cell is attached to the $n-1$ skeleton $K^{[n-1]}$ by an
attaching map
$$ a\co  S^{n-1} \rightarrow K^{[n-1]}. $$
The following definition was the subject of \cite[Section~4]{Behrensroot}.

\begin{defn}
Let $\beta$ be an element of $\pi_j(S)$.  Then the $K$--Toda bracket
$\bra{K}(\beta)$ is a lift
$$ \xymatrix{
S^{j+n-1} \ar[r]^\beta \ar@{.>}[rrd]_{\bra{K}(\beta)} &
S^{n-1} \ar[r]^a &
K^{[n-1]}
\\
&& S^0 \ar[u]_\iota
} $$
The $K$--Toda bracket may not exist, or may not be well defined.
\end{defn}

Given a $k_i^\mathrm{th}$ filtered root invariant, 
the $k_{i+1}^\mathrm{th}$ filtered root
invariant may be revealed by the presence of an Adams differential.

\begin{thm}{\rm\cite[Theorem~5.3]{Behrensroot}}\label{thmB}\qua
There is a (possibly trivial) $E$--Adams differential
$$ d_r(R^{[k_i]}_{E}(\alpha)) =
  \bigl\langle{P_{-N_i}^{-N_{i+1}}}\bigr\rangle
  \bigl(R^{[k_{i+1}]}_E(\alpha)\bigr).$$
\end{thm}

\fullref{thmB} is saying the following differential 
happens in the $E$--ASS
chart.
$$
\xymatrix{
&& \bullet 
\\
R^{[k_{i+1}]}_E(\alpha)
\ar@{-}[rru]^{\bigl\langle{P_{-N_i}^{-N_{i+1}}}\bigr\rangle} 
\\
&&& R^{[k_i]}_E(\alpha) \ar[uul]
}
$$
If this differential is zero, there may still be a hidden extension that
reveals the $k_{i+1}^\mathrm{th}$ filtered root invariant.  In the next
theorem, we use the notation $\wbar{\beta}$ for an element that $\beta \in
E_1^{*,*}$ detects in the $E$--ASS.

\begin{thm}{\rm\cite[Theorem~5.4]{Behrensroot}}\label{thmC}\qua
There is an equality of (possibly trivial) elements of $\pi_*(S)$
$$\bigl\langle{P^{-N_i}_{-M}}\bigr\rangle
  \bigl(\br{R^{[k_i]}_E(\alpha)}\bigr) = 
  \br{\bigl\langle{P_{-M}^{-N_{i+1}}}\bigr\rangle
  \bigl(R^{[k_{i+1}]}_E(\alpha)\bigr)}. $$
Here $M$ is the largest integer for which 
$\bigl\langle{P^{-N_i}_{-M}}\bigr\rangle
 \bigl(\br{R^{[k_i]}_E(\alpha)}\bigr)$ exists and is non-trivial,
and the second Toda bracket is taken in the $E$--ASS.
\end{thm}

\fullref{thmC} says that the following hidden extension
happens on the $E$--ASS chart.
$$\xymatrix{
&&& \bullet
\\
R^{[k_{i+1}]}_E(\alpha)
  \ar@{-}[rrru]^{\bigl\langle{P^{-N_{i+1}}_{-M}}\bigr\rangle}
\\
& R^{[k_{i}]}_E(\alpha)
  \ar@{--}[uurr]_{\bigl\langle{P^{-N_i}_{-M}}\bigr\rangle}
}
$$

We have given tools to move from one filtered root invariant to the next,
but we need a place to start this process.  For $E = BP$ in
\cite{Behrensroot} we used $BP$--root invariants and $BP \wedge
\br{BP}$--root invariants.  For $E = H\FF_p$ the first filtered root
invariant is given by the algebraic root invariant $R_{\alg}$
(see, for example Mahowald--Ravenel \cite{MahowaldRavenel}).  The nice
thing about $R_{\alg}$
is that it is very computable, especially with the help of a computer 
(see Bruner \cite{Bruner}). 

\begin{thm}{\rm\cite[Theorem~5.10]{Behrensroot}}\label{thmD}\qua
Let $\alpha$ be of Adams filtration $s$, detected in the ASS 
by $\wtilde{\alpha}$ in
$\ext$.  Then
the first $H\FF_p$--filtered root invariant is given by
$$ R^{[s]}_{H\FF_p}(\alpha) = R_{\alg}(\alpha). $$
\end{thm}

\part{$2$--primary calculations}
\setobjecttype{Part}
\label{part:calc}

In this part we are always implicitly working $2$--locally.  Our goal is to
explain how the theory of \fullref{part:lec} play out in low dimensions in
the ANSS and the ASS at the prime $2$.  Unlike in \fullref{part:lec}, we
intend to be completely precise about these calculations.  This part is
really an extension of \cite{Behrensroot} to the prime $2$.

Our main result is to compute the
homotopy root invariants of all of the $v_1$--periodic elements through the
$12$--stem (\fullref{thm:hroots}).  These root invariants turn out fit
into the primary $v_2$--family investigated by Mahowald \cite{Mahowaldv2}
(but beware! $v_2^8$ does not exist, $v_2^{32}$ does exist 
; see Hopkins--Mahowald \cite{HopkinsMahowaldself}).

Our plan of attack is the following.  In \cite{Behrensroot}, we computed the
$BP\wedge \td{BP}$--root invariants of the elements $\alpha_{i/j}$.  These
give the first filtered root invariant modulo an indeterminacy group 
as described
in \fullref{prop:BP^1root}.  
Once we know the indeterminacy group, we can
identify $R^{[1]}_{BP}(\alpha_{i/j})$, and then get
$R^{[2]}_{BP}(\alpha_{i/j})$.  The higher filtered root invariants are
deduced from differentials and hidden extensions in the ANSS.  We then
check to see that these top filtered root invariants must detect the
homotopy root invariants.

In \fullref{sec:indeterminacy}, we compute this indeterminacy group
completely.  This essentially involves understanding how the $v_1$--periodic
elements act in the algebraic Atiyah--Hirzebruch spectral sequence for
$\ext(BP_*P^\infty)$.  Unfortunately, there is no $J$--homomorphism in
$\ext$ to
produce these differentials, so we must resort to explicit computation of
the differentials using the $BP_*BP$--coaction on $BP_*P^\infty$.  We find
generators for this $BP$--homology group that make a complete determination
of the AAHSS differentials possible.  This method may be interesting in its
own right, in the sense that it gives a particularly clean and 
pleasant description of the coaction.

In \fullref{sec:v0filt} we describe how the theory of $BP$--filtered
root invariants reproduces the expected root invariants of the elements 
$2^i$.  The
differentials and hidden extensions amongst the elements
$\alpha_{i/j}\alpha_1^k$ get a rather natural interpretation: the ANSS must
deal with the fact that the homotopy Greek letter elements $\alpha_i^h$ are
different from the algebraic Greek letter elements $\alpha_i$.

In \fullref{sec:v1filt}, we use our computation of the indeterminacy
group to compute the the filtered root invariants $R^{[k]}_{BP}(\alpha_{i/j})$
for $k = 1,2$.  We find that the indeterminacy is essential to allow for
the root invariants $R(\eta) = \nu$ and $R(\nu) = \sigma$.  For the higher
dimensional $v_1$--periodic elements, we find that for all other $i$ and $j$,
$$ R^{[2]}_{BP}(\alpha_{i/j}) \doteq \beta_{i/j} + \text{something} $$
where the ``something'' is rather innocuous.  One could take this
calculation as further evidence that the root invariants wants to take
$v_n$--periodic families of Greek letter elements to $v_{n+1}$--periodic
families of
Greek letter elements.

In \fullref{sec:froots} we compute all of the higher $BP$--filtered root
invariants of the elements $\alpha_{i/j}$ which lie within the $12$--stem.
These higher filtered root invariants are deduced from differentials and
hidden extensions in the ANSS, and, amusingly enough, actually account for
most of the differentials and hidden extensions in this range.  We saw this
sort of behavior at the prime $3$ in \cite{Behrensroot}.  The reason the range
is so limited is the author's limited knowledge of the ANSS at the prime
$2$.

In \fullref{sec:hroots} we show that the filtered root invariants of
\fullref{sec:froots} actually detect homotopy root invariants.  This is
done by brute force.
We show that there are no elements of
$\pi_*(P_{-N}^\infty)$ that could survive to the
difference of the filtered root invariant and the homotopy root invariant.
We use the $BP$--filtered root invariants for the 
elements in $BP$--Adams filtration $1$, and the $H\FF_2$--filtered root
invariants for the rest.

\section{The indeterminacy spectral sequence}\label{sec:indeterminacy}

Recall that
$$ BP_*(P^{2k}_{2l-1}) = BP_*\{ e_{2m-1} \: : \: l \le m \le k \}/\Bigl(\sum_{i
\ge 0} c_ie_{2(m-i)-1}\Bigr)$$
where the universal $2$--typical $2$--series is given by
$$ [2]_F(x) = \sum_{i \ge 0}c_i x^{i+1} \in BP_*\llbracket x\rrbracket. $$
In particular, the first couple of values of $c_i$ are 
$$c_0 = 2, \qquad c_1 = -v_1.$$
We will first define a $v_1$--self map of $BP_*BP$--comodules.  Define a map
$$ \wtilde{v}_1 \co  BP_*P_{2l+1}^{2k} \rightarrow BP_*P_{2l-1}^{2(k-1)} $$
by
$$ \wtilde{v}_1(e_{2m-1}) = \sum_{i \ge 1} -c_i e_{2(m-i)-1} $$
This is a map of comodules since in
$BP_*P_{2l-1}^{2k}$, we have
$$ \wtilde{v}_1(e_{2m-1}) = 2e_{2m-1}. $$
Thus $\wtilde{v}_1$ is just a certain factorization of multiplication by $2$.

The short exact sequences
$$ 0 \rightarrow BP_*P^{2(k-1)} \rightarrow BP_*P^{2k} \rightarrow
\Sigma^{2k-1} BP_*/(2) \rightarrow 0 $$
gives rise to long exact sequences of Ext groups, which piece together to
give a \emph{modified AAHSS}
$$ E_1^{k,m,s} = \ext^{s,s+k}(\Sigma^{2m-1}BP_*/(2)) \Rightarrow
\ext^{s,s+k}(BP_*P^\infty). $$
We shall refer to elements of $E_1^{k,m,s}$ of the modified AAHSS by
$x[2m-1]$, where $x$ is an element of
$\ext^{s,k+s}(\Sigma^{2m-1}BP_*/(2))$.
The existence of the map $\wtilde{v}_1$ gives the following propagation result in
the modified AAHSS.

\begin{prop}\label{prop:v1prop}
Suppose $x[2m-1]$ is an element of the modified AAHSS, and that there
is a differential
$$ d_r(x[2m-1]) = y[2(m-r)-1]. $$
Then we have
$$ d_r(v_1x[2(m-1)-1]) = v_1y[2(m-1-r)-1]. $$
\end{prop}

\begin{proof}
The map $\wtilde{v}_1$ induces a map of modified AAHSS's:
$$
\xymatrix{
\bigoplus_{l+1 \le m \le k} \ext^{s,s+k}(\Sigma^{2m-1}BP_*/(2))
\ar@{=>}[r] \ar[d]_{v_1} &
\ext^{s,s+k}\bigl(BP_*P_{2l+1}^{2k}\bigr) \ar[d]^{(\wtilde{v}_1)_*}
\\
\bigoplus_{l+1 \le m \le k} \ext^{s,s+k}(\Sigma^{2m-3}BP_*/(2)) 
\ar@{=>}[r] & 
\ext^{s,s+k}\bigl(BP_*P_{2l-1}^{2(k-1)}\bigr)
}$$
This proves the proposition.
\end{proof}

\begin{prop}\label{prop:AAHSSgendiffs}
The differentials on the elements $1[2m-1]$ in the modified AAHSS are given
as follows:
\begin{alignat*}{2}
d_1(1[2m-1]) & = \alpha_1[2(m-1)-1] & \qquad & \text{$m$ odd} \\
d_2(1[2m-1]) & = \wtilde{\beta}_1[2(m-2)-1] & \qquad & \nu_2(m) = 1 \\
d_r(1[2m-1]) & = v_1^{k-2} (x_7+\wtilde{\beta}_{2/2})[2(m-k-2)-1] & 
\qquad & \nu_2(m) = k, k = 2,3 \\
d_r(1[2m-1]) & = v_1^{k-2} x_7[2(m-k-2)-1] & 
\qquad & \nu_2(m) = k, k \ge 4 \\
\end{alignat*}
\end{prop}

\begin{proof}
The formulas for $d_1$ and $d_2$ follow immediately from the well known
attaching map structure of $P^\infty$.
We shall prove the formulas for the higher differentials 
by working with the negative cells of
$P_{-N}^\infty$, and then by using James periodicity.  It suffices to
consider $m = -2^k$.
There is the following equivalence to the Spanier--Whitehead dual
(see Bruner--May--McClure--Steinberger \cite{Hinfty}).
$$ \Sigma P^{-2\cdot 2^k}_{-2(2^k+l)-1} \simeq DP_{2\cdot 2^k-1}^{2(2^k+l)} $$
It follows that we have an isomorphism of $BP_*BP$--comodules
$$ BP_*(\Sigma P^{-2\cdot 2^k}_{-2(2^k+l)-1}) \cong BP^{-*}P_{2\cdot
2^k-1}^{2(2^k+l)}. $$
Here, for finite $X$, the cohomology group $BP^{-*}X$ is viewed as a 
$BP_*BP$--comodule by the coaction given by the
composite
\begin{align*}
BP^{-*}X = \pi_*(F(X,BP)) \xrightarrow{(\eta_R)_*} \: & \pi_*(F(X,BP \wedge
BP)) \\
& \cong \pi_*(BP \wedge BP \wedge DX) \\
& \cong BP_*BP \otimes_{BP_*} BP^{-*}X.
\end{align*}
We recall from \cite{Behrensroot} that there are short exact sequences
$$ 0 \rightarrow 
BP^*\CC P_a^b \xrightarrow{\cdot[2]_F(x)/x} BP^*\CC P_a^b \rightarrow
BP^*P_{2a-1}^{2b} \rightarrow 0. $$
Here we have
$$ BP^{-*}\CC P^\infty \cong BP_*\llbracket x \rrbracket $$
where $x$ has (homological degree) $-2$, and $BP^{-*}\CC P_a^b$ is given by
the ideal
$$ BP^{-*}\CC P_a^b \cong (x^a) \subseteq BP_*[x]/(x^{b+1}) 
\cong BP^{-*}\CC P^b. $$
We recall from \cite{Behrensroot} that the coaction of $BP_*BP$ on $h(x) \in
BP^{-*} \CC P_a^b $ is given by
$$ \psi(h(x)) = (f_*h)(f(x)) $$
where $f$ is the universal isomorphism of $2$--typical formal groups, whose
inverse is given by
$$ f^{-1}(x) = \sum_{i \ge 0}^F t_i x^{2^i}. $$
The polynomial $f_*h(x)$ is the polynomial obtained by applying the right
unit to all of the coefficients of $h(x)$.  The surjection of $BP^*\CC
P_a^b$ onto $BP^*P_{2a-1}^{2b}$ completely determines the latter as a
$BP_*BP$--comodule.  In what follows we shall refer to elements of
$BP^*P_{2a-1}^{2b}$ by the names of elements in $BP^*\CC P_a^b$ which project
onto them.

We shall need the following formulas. (We are using Hazewinkel generators.)
\begin{align*}
f(x) & = x - t_1 x^2 + (2t_1^2+v_1t_1)x^3 + \cdots \\
[2]_F(x) & = 2x - v_1x^2 + 2v_1^2x^3 + \cdots \\
\eta_R(v_1) & = v_1 + 2t_1 \\
\eta_R(v_2) & = v_2 + 2t_2 - 4t_1^3 - 5v_1t_1^2 - 3v_1^2t_1
\end{align*}
We shall now use our very specific knowledge of the $BP_*BP$ coaction to
determine the differential $d(1[-2\cdot 2^k-1])$ in the modified AAHSS for
$$ \ext\bigl(BP_*\Sigma P_{-2(2^k-l)-1}^{-2\cdot 2^k}\bigr) = 
\ext\bigl(BP^{-*}P^{2(2^k+l)}_{2\cdot 2^k-1}\bigr). $$
We do this for $l = k+2$, so in what follows 
we work modulo $x^{2^k+k+3}$.  The desired
differential is governed by the coaction on $x^{2^k} \in 
\smash{BP^{-*}P^{2(2^k+l)}_{2\cdot 2^k-1}}$.  Actually, the case
$k=2$ must be handled separately, because in the computations that follow
we are implicitly using the fact that $2^k > k+2$.  However, the method,
and conclusion,
for $k=2$ are completely identical.
\begin{eqnarray*}
\psi\bigl(x^{2^k}\bigr) & = & \bigl(f(x)\bigr)^{2^k} \\
& = & \bigl(x - t_1 x^2 + \bigl(2t_1^2+v_1t_1\bigr)x^3 + \cdots\bigr)^{2^k} \\
& = & x^{2^k} - 
\tbinom{2^k}{1}t_1x^{2^k+1} + 
\tbinom{2^k}{1}\bigl(2t_1^2+v_1t_1\bigr)x^{2^k+2} + \\
&& \tbinom{2^k}{2}t_1^2x^{2^k+2} + 
\tbinom{2^k}{4}t_1^4x^{2^k+4} +
\cdots \\
& = & x^{2^k} - 
2^kt_1x^{2^k+1} + 
2^{k+1}t_1^2x^{2^k+2} + 
2^kv_1t_1x^{2^k+2} + \\
&& \bigl(2^k-1\bigr)2^{k-1}t_1^2x^{2^k+2} +
 2^{k-2}\cdot\tfrac13\bigl(2^k-1\bigr)\bigl(2^{k-1}-1\bigr)\bigl(2^k-3\bigr)
  t_1^4 x^{2^k+4} +
\cdots \\
& = & x^{2^k} - 
v_1^kt_1x^{2^k+k+1} - 
v_1^{k-1}t_1^2x^{2^k+k+1} +
v_1^{k-2}t_1^4 x^{2^k+k+2} + \\
&& v_1^{k+1}t_1x^{2^k+k+2} + 
\cdots
\end{eqnarray*}
We conclude that in the cobar complex for
$\smash{\ext\bigl(BP^*P_{2\cdot2^k-1}^{2(2^k+k+2)}\bigr)}$, we have:
$$ d\bigl(x^{2^k}\bigr) {=} v_1^kt_1x^{2^k+k+1} {+} 
v_1^{k-1}t_1^2x^{2^k+k+1} {+}
v_1^{k-2}t_1^4 x^{2^k+k+2} {+}
v_1^{k+1}t_1x^{2^k+k+2} {+} \cdots$$
We compute the differential in the modified AAHSS by adding a coboundary
supported on an element of lower cellular filtration.  Namely, we compute
the coaction on $v_1^{k-2}v_2x^{2^k+k+1}$ as:
\begin{align*}
\psi\bigl(v_1^{k-2} & v_2x^{2^k+k+1}\bigr) = 
\eta_R\bigl(v_1^{k-2}v_2\bigr)\bigl(f(x)\bigr)^{2^k+k+1} 
\\
= &
(v_1+2t_1)^{k-2}\bigl(v_2 + 2t_2 - 4t_1^3 - 5v_1t_1^2 - 3v_1^2t_1\bigr)
\bigl(x-t_1x^2+\cdots\bigr)^{2^k+k+1}
\\
= & 
v_1^{k-2}\bigl(v_2 + 2t_2 - 4t_1^3 - 5v_1t_1^2 - 3v_1^2t_1\bigr)x^{2^k+k+1} + \\
& 2(k-2)v_1^{k-3}t_1\bigl(v_2 + 2t_2 - 4t_1^3 - 5v_1t_1^2 - 
 3v_1^2t_1\bigr)x^{2^k+k+1} - \\
& (2^k+k+1)v_1^{k-2}t_1\bigl(v_2 + 2t_2 - 4t_1^3 - 5v_1t_1^2 - 
 3v_1^2t_1\bigr)x^{2^k+k+2} 
+ \cdots
\\
= & 
\bigl(v_1^{k-2}v_2 + v_1^{k-1}t_1^2 + v_1^kt_1\bigr)x^{2^k+k+1} +
v_1^{k-1}\bigl(t_2 + v_1t_1^2\bigr)x^{2^k+k+2} + \\
& \bigl(2^k+2k-1\bigr)v_1^{k-2}t_1\bigl(v_2 + 2t_2 - 4t_1^3 - 5v_1t_1^2 - 
 3v_1^2t_1\bigr)x^{2^k+k+2} 
+ \cdots \\
= & 
\bigl(v_1^{k-2}v_2 + v_1^{k-1}t_1^2 + v_1^kt_1\bigr)x^{2^k+k+1} + \\
& \bigl(v_1^{k-1}t_2 + v_1^{k-2}v_2t_1 + v_1^{k-1}t_1^3\bigr)x^{2^k+k+2} 
+ \cdots
\end{align*}
We conclude that in the cobar complex for $\smash{\ext\bigl(BP^*P_{2\cdot
2^k-1}^{2(2^k+k+2)}\bigr)}$ we have:
\begin{multline*} 
d\bigl(v_1^{k-2}v_2x^{2^k+k+1}\bigr) = \\
\bigl(v_1^{k-1}t_1^2 + v_1^kt_1\bigr)x^{2^k+k+1} +
\bigl(v_1^{k-1}t_2 + v_1^{k-2}v_2t_1 + v_1^{k-1}t_1^3\bigr)x^{2^k+k+2} + \cdots
\end{multline*}
We therefore have:
\begin{multline*}
d\bigl(x^{2^k} - v_1^{k-2}v_2x^{2^k+k+1}\bigr) = \\
\bigl(v_1^{k-2}t_1^4 + v_1^{k+1}t_1 + v_1^{k-1}t_2 + 
v_1^{k-2}v_2t_1 + v_1^{k-1}t_1^3\bigr)x^{2^k+k+2}
\end{multline*}
We recall from Ravenel \cite{Ravenel} that the generators of $\ext^{1,8}(BP_*/(2))$
are $x_7$ and $\wtilde{\beta}_{2/2}$, and they are represented 
in the cobar complex by the elements
\begin{align*}
x_7 & = v_1t_2 + v_2t_1 + v_1t_1^3 \\
\wtilde{\beta}_{2/2} & = t_1^4 + v_1^3t_1
\end{align*}
We conclude that for $k = 2,3$ we have the modified AAHSS differentials
$$ d_{k+2}\bigl(1[-2\cdot2^k-1]\bigr) = v_1^{k-2}\bigl(x_7 +
\wtilde{\beta}_{2/2}\bigr)\bigl[-2\bigl(2^k+k+2\bigr)-1\bigr]. $$
If $k \ge 4$, then we may add an additional 
coboundary to obtain the cobar formula
$$ d\bigl(x^{2^k} + v_1^{k-2}v_2x^{2^k+k+1} +
v_1^{k-4}v_2^2x^{2^k+k+2}\bigr) =
x_7 x^{2^k+k+2} $$
from which it follows that for $k \ge 4$, we have the modified AAHSS
differential
$$ d_{k+2}\bigl(1\bigl[-2\cdot2^k-1\bigr]\bigr) =
v_1^{k-2}x_7\bigl[-2\bigl(2^k+k+2\bigr)-1\bigr].\proved$$
\end{proof}

Combining \fullref{prop:AAHSSgendiffs} with
\fullref{prop:v1prop}, we get the following differentials.

\begin{prop}\label{prop:AAHSSv1diffs}
The differentials on the elements $v_1^i[2m-1]$ for $i \ge 1$
in the modified AAHSS are given
as follows.
\begin{alignat*}{2}
d_1\bigl(v_1^i[2m-1]\bigr) & = v_1^i\alpha_1[2(m-1)-1] & \qquad & \text{$m+i$ odd} \\
d_3\bigl(v_1^i[2m-1]\bigr) & = 
v_1^{i-1} x_7[2(m-3)-1] & \qquad & \nu_2(m+i) = 1 \\
d_4\bigl(v_1^i[2m-1]\bigr) & = 
v_1^{i} (x_7+\beta_{2/2})[2(m-4)-1] & \qquad & \nu_2(m+i) = 2 \\
d_r\bigl(v_1^i[2m-1]\bigr) & = 
v_1^{k+i-2}x_7[2(m-k-2)-1] & \qquad & \nu_2(m+i) = k, k \ge 3
\end{alignat*}
\end{prop}

\begin{proof}
The differentials follow from applying
$v_1$--propagation as described in \fullref{prop:v1prop} to the
differentials of \fullref{prop:AAHSSgendiffs}.  However,
the element 
$v_1\wtilde{\beta}_1$ is null in $\ext(BP_*/(2))$.  An explicit computation
similar to that in the proof of \fullref{prop:AAHSSgendiffs} yields
the modified AAHSS differential
$$ d_3(v_1[2m-1]) = x_7[2(m-3)-1] $$
for $m+1 \equiv 2 \pmod 4$.  The rest of the $d_3$'s then follow by
$v_1$--propagation.
\end{proof}

Recall that we have the following computation, which is given by
combining Corollary~$5.9$ and Proposition~$10.2$ of \cite{Behrensroot}.

\begin{prop}\label{prop:BP^1root}
There is an indeterminacy group $A_{i/j} \subseteq BP_*\td{BP}$ such that
$$ R^{[1]}_{BP}(\alpha_{i/j}) \subseteq \wtilde{\beta}_{i/j} + 2BP_*\td{BP} +
A_{i/j}. $$
\end{prop}

In \cite{Behrensroot} (spectral sequence~(10.8) and the discussion which
follows it), a method of computing
this indeterminacy group was described in terms of the differentials of 
an \emph{indeterminacy spectral sequence}.  
This spectral sequence is a truncated
version of the AAHSS.  The differentials given in
\fullref{prop:AAHSSgendiffs} and
\fullref{prop:AAHSSv1diffs} give differentials in the indeterminacy
spectral sequence, which translate to the following result.

\begin{prop}\label{prop:indeterminacy}
The indeterminacy group $A_{i/j}$ for 
$R_{BP}^{[1]}(\alpha_{i/j})$ is contained in the $\ZZ_{(2)}$
module spanned by $2BP_*BP$ and the generator given in the table below,
where $a = 3i-j$.

\begin{center}
\begin{tabular}{|c|c|l|}
\hline
$a$ & Generator & Condition \\
\hline \hline
$2$ & $-$ & $\nu_2(i) \ge 2$ \\
    & $\wtilde{\beta}_1$ & $\nu_2(i) = 1$ \\
    & $v_1\alpha_1$ & $\nu_2(i) = 0$ \\
\hline
$3$ & $-$ & $\nu_2(i) \ge 1$ \\
    & $v_1^2\alpha_1$ & $\nu_2(i) = 0$ \\
\hline
$4$ & $-$ & $\nu_2(i) \ge 3$ \\
    & $x_7+\wtilde{\beta}_{2/2}$ & $\nu_2(i) = 2$ \\
    & $x_7$ & $\nu_2(i) = 1$ \\
    & $v_1^3\alpha_1$ & $\nu_2(i) = 0$ \\
\hline
$5$ & $-$ & $\nu_2(i) \ge 4$ \\
    & $v_1(x_7+\wtilde{\beta}_{2/2})$ & $2 \le \nu_2(i) \le 3$ \\
    & $v_1x_7$ & $\nu_2(i) = 1$ \\
    & $v_1^4\alpha_1$ & $\nu_2(i) = 0$ \\
\hline
$a \ge 6$ & $-$ & $\nu_2(i) \ge a-1$ \\
    & $v_1^{a-4}x_7$ & $1 \le \nu_2(i) \le a-2$ \\
    & $v_1^{a-1}\alpha_1$ & $\nu_2(i) = 0$ \\
\hline
\end{tabular}
\end{center}

\end{prop}

\begin{rmk}
Not all of the entries of the table in \fullref{prop:indeterminacy}
actually occur with the allowable values of $i$ and $j$ for $\alpha_{i/j}$.
\end{rmk}

\section{$BP$ filtered root invariants of $2^k$}\label{sec:v0filt}

The root invariants of the elements $2^k$ were determined by
Mahowald and Ravenel \cite{MahowaldRavenel} and by Johnson \cite{Johnson}.
In this section we will explain how the root invariants of the
elements $2^k$ are formed from the perspective of the ANSS.  This analysis
provides an explanation for the pattern of differentials and hidden
extensions in the $v_1$--periodic part of the ANSS.  We only compute
filtered root invariants --- 
our treatment is not an independent
verification of the know values of the homotopy root invariants 
$R(2^k)$ because we do not eliminate the
higher Adams--Novikov filtration obstructions required by
\fullref{thmA}.  The proper thing to do would be to combine the use of
the ASS and the ANSS, in a manner employed in the infinite family
computations of \cite{Behrensroot}.

\begin{figure}[htp!]
\labellist
\tiny
\pinlabel {$0$} [tr] at 9 11

\pinlabel {$1$} [r] at 9 50
\pinlabel {$2$} [r] at 9 86
\pinlabel {$3$} [r] at 9 122
\pinlabel {$4$} [r] at 9 158
\pinlabel {$5$} [r] at 9 194
\pinlabel {$6$} [r] at 9 230
\pinlabel {$7$} [r] at 9 266
\pinlabel {$8$} [r] at 9 302

\pinlabel  {$1$} [t] <1pt,0pt> at 30 11
\pinlabel  {$2$} [t] <1pt,0pt> at 48 11
\pinlabel  {$3$} [t] <1pt,0pt> at 66 11
\pinlabel  {$4$} [t] <1pt,0pt> at 84 11
\pinlabel  {$5$} [t] <1pt,0pt> at 102 11
\pinlabel  {$6$} [t] <1pt,0pt> at 120 11
\pinlabel  {$7$} [t] <1pt,0pt> at 138 11
\pinlabel  {$8$} [t] <1pt,0pt> at 156 11
\pinlabel  {$9$} [t] <1pt,0pt> at 174 11
\pinlabel {$10$} [t] <1pt,0pt> at 192 11
\pinlabel {$11$} [t] <1pt,0pt> at 210 11
\pinlabel {$12$} [t] <1pt,0pt> at 228 11
\pinlabel {$13$} [t] <1pt,0pt> at 246 11
\pinlabel {$14$} [t] <1pt,0pt> at 264 11
\pinlabel {$15$} [t] <1pt,0pt> at 282 11
\pinlabel {$16$} [t] <1pt,0pt> at 300 11
\pinlabel {$17$} [t] <1pt,0pt> at 318 11
\pinlabel {$18$} [t] <1pt,0pt> at 336 11
\pinlabel {$19$} [t] <1pt,0pt> at 354 11
\pinlabel {$20$} [t] <1pt,0pt> at 372 11
\pinlabel {$21$} [t] <1pt,0pt> at 390 11
\pinlabel {$22$} [t] <1pt,0pt> at 408 11
\pinlabel {$23$} [t] <1pt,0pt> at 426 11
\pinlabel {$24$} [t] <1pt,0pt> at 444 11
\pinlabel {$25$} [t] <1pt,0pt> at 462 11
\pinlabel {$26$} [t] <1pt,0pt> at 480 11
\pinlabel {$27$} [t] <1pt,0pt> at 498 11
\pinlabel {$28$} [t] <1pt,0pt> at 516 11
\pinlabel {$29$} [t] <1pt,0pt> at 534 11
\pinlabel {$30$} [t] <1pt,0pt> at 552 11

\pinlabel {$\infty$} at 14 14
\pinlabel {$2$} at 67 50
\pinlabel {$4$} at 139 50
\pinlabel {$3$} at 211 50
\pinlabel {$5$} at 283 50
\pinlabel {$3$} at 355 50
\pinlabel {$4$} at 427 50
\pinlabel {$3$} at 500 50
\pinlabel {$2$} at 337 77
\pinlabel {$2$} at 373 163
\pinlabel {$3$} at 427 194

\pinlabel {$\alpha_1$} [l] at 31 50
\pinlabel {$\alpha_2$} [l] at 72 50
\pinlabel {$\alpha_3$} [l] at 103 50
\pinlabel {$\alpha_4$} [l] at 144 50
\pinlabel {$\alpha_5$} [l] at 175 50
\pinlabel {$\alpha_6$} [l] at 216 50
\pinlabel {$\alpha_7$} [l] at 247 50
\pinlabel {$\alpha_8$} [l] at 289 50
\pinlabel {$\alpha_9$} [l] at 319 50
\pinlabel {$\alpha_{10}$} [l] at 361 50
\pinlabel {$\alpha_{11}$} [l] at 391 50
\pinlabel {$\alpha_{12}$} [l] at 433 50
\pinlabel {$\alpha_{13}$} [l] at 463 50
\pinlabel {$\alpha_{14}$} [l] at 505 50
\pinlabel {$\alpha_{15}$} [l] at 535 50

\pinlabel {$\beta_{2\!/\!2}$} [b] at 121 86
\pinlabel {$\beta_2$} [br] at 158 87
\pinlabel {$\beta_{4\!/\!4}$} [t] at 266 87
\pinlabel {$\beta_3$} [br] at 266 90
\pinlabel {$\beta_{4\!/\!3}$} [tl] at 300 83
\pinlabel {$\beta_{4\!/\!2}$} [t] at 338 74
\pinlabel {$\beta_4$} [r] at 374 86
\pinlabel {$\eta_2$} [tl] at 360 117
\pinlabel {$\eta_{3\!/\!2}$} [t] at 428 123
\pinlabel {$\beta_5$} [r] at 482 92
\pinlabel {$\beta_{8\!/\!8}$} [tl] at 554 86
\pinlabel {$\beta_{6\!/\!2}$} [bl] at 554 92
\pinlabel {$\Pi\beta_{2\!/\!2}$} [br] at 482 230
\pinlabel {$\Pi\beta_2$} [t] at 518 231
\endlabellist
\centerline{\rotatebox{90}{\includegraphics[width=7.10in]{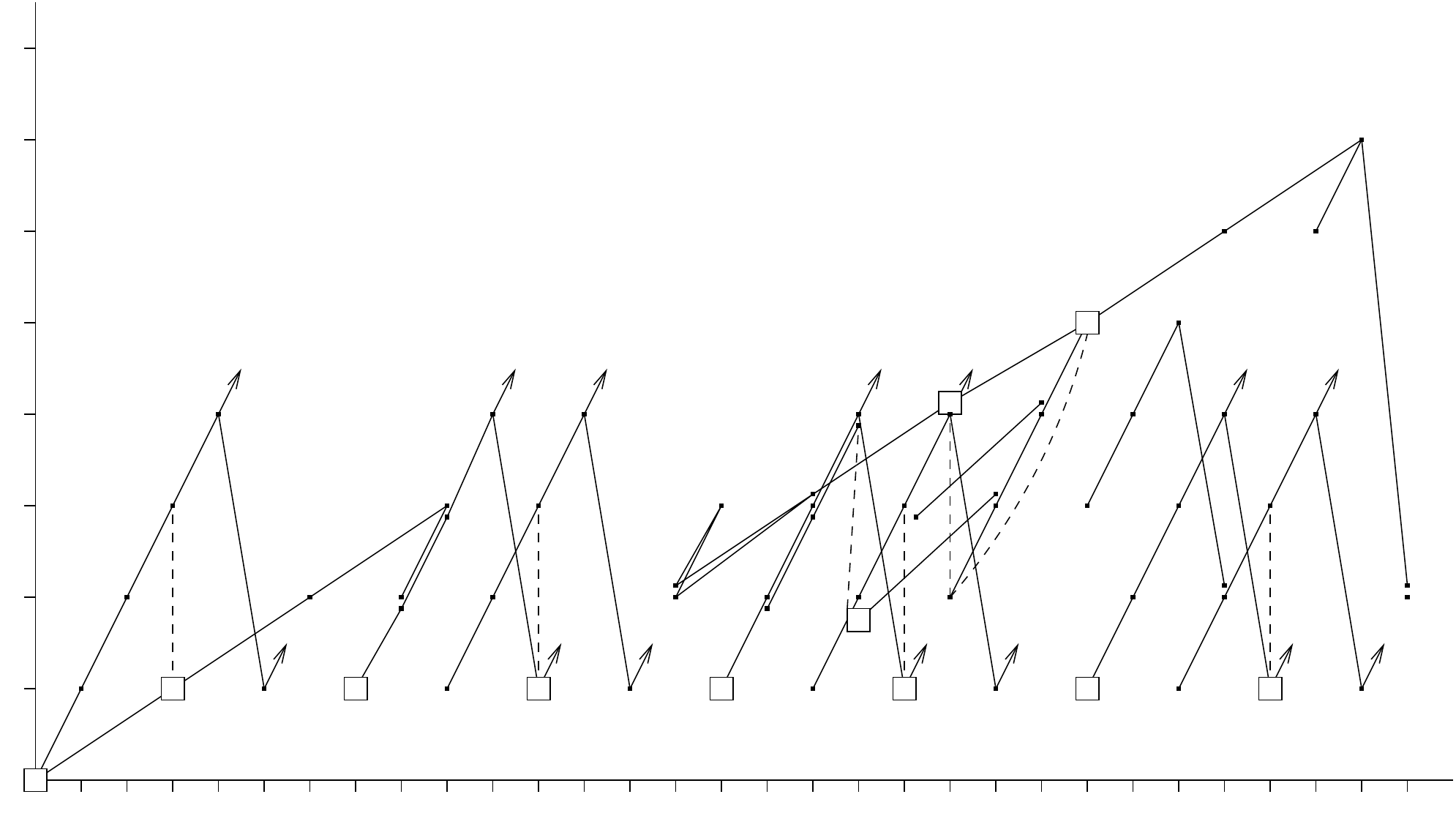}}}
\caption{The Adams--Novikov spectral sequence at $p=2$}
\label{fig:anss2}
\end{figure}
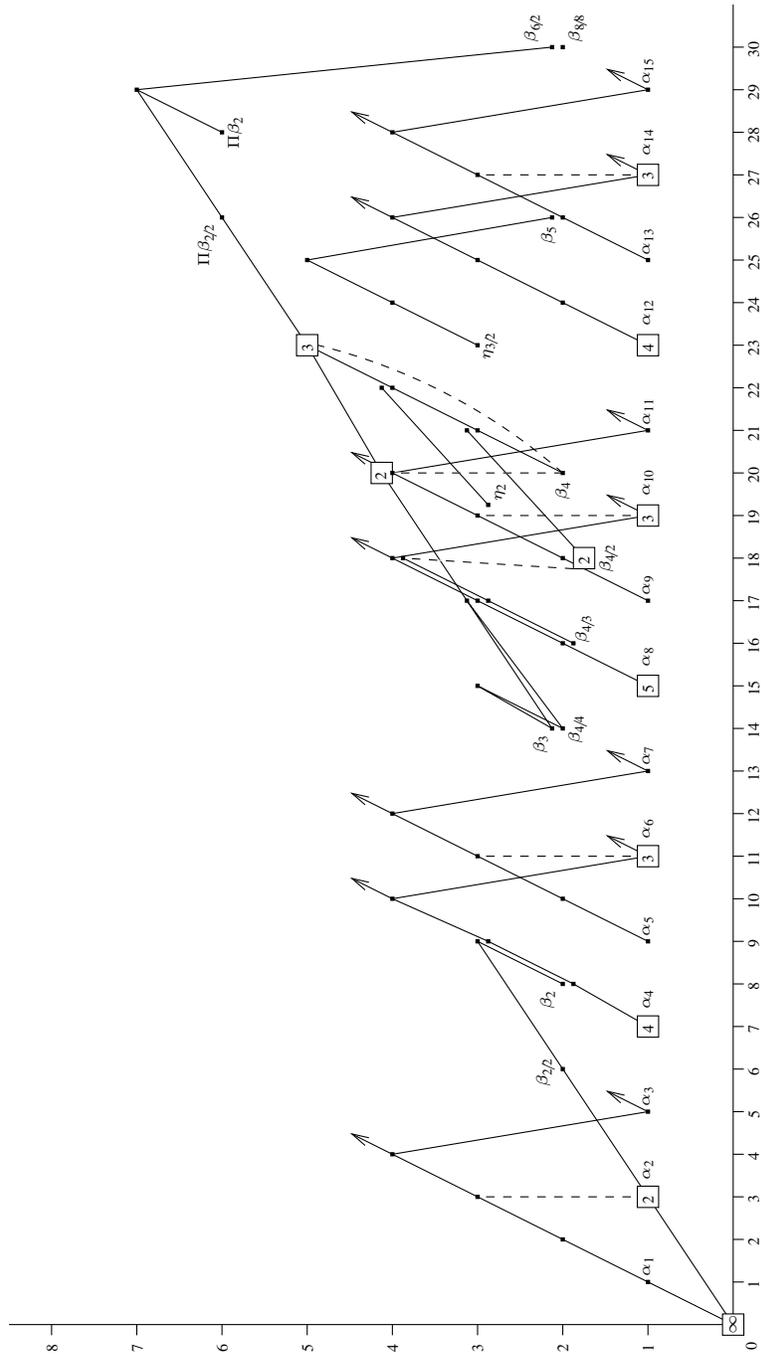

Throughout this section and the rest of the paper, the reader should find 
it helpful to refer to 
\fullref{fig:anss2}, which depicts the ANSS chart for the sphere at the
prime $2$ through the $29$--stem.

\cite[Proposition~10.1]{Behrensroot} states that 
$$ \alpha_k \in R^{[1]}_{BP}(2^k). $$
Our description of $R(2^k)$ depends on the value of $k$ modulo $4$. 
\begin{description}
\item[${k \equiv 1 \pmod 4}$] The element $\alpha_k$ is a
permanent cycle which detects the homotopy root invariant $R(2^k)$.

\item[${k \equiv 2 \pmod 4}$] While $\alpha_k$ is a permanent cycle,
it does not detect an element of $R(2^k)$.  This is an instance
where \fullref{thmC} comes into play.  
The $-2k$--cell attaches to the $(-2k{-}1)$--cell of $P^{-2k}_{-2k-1}$ 
by the degree $2$
map, and the
$(-2k{+}1)$--cell attaches to the $(-2k{-}1)$--cell in $P^{-2k+1}_{-2k-1}$ by the map 
$\alpha_1$.
There is a hidden
extension $\wbar{\alpha}_{1} \cdot \wbar{\alpha}_{k-1}\wbar{\alpha}_1 = 2 \cdot 
\wbar{\alpha}_k$, so we may use \fullref{thmC} to deduce that
$$ \alpha_{k-1} \alpha_1 \in R^{[2]}_{BP}\bigl(2^k\bigr). $$
The element $\alpha_{k-1}\alpha_1$ detects the homotopy root invariant
$R\bigl(2^k\bigr)$.

\item[${k \equiv 3 \pmod 4}$] 
The first filtered root invariant $\alpha_k$ is not a permanent cycle, 
so we turn to 
\fullref{thmB}.  There is an Adams--Novikov differential
$$ d_3(\alpha_k) = \alpha_1 \cdot \alpha_{k-2} \alpha_1^2. $$
The $(-2k{+}2)$--cell attaches to the $-2k$--cell in $P^{-2k+2}_{-2k}$ 
with attaching map $\alpha_1$.
We conclude that 
$$ \alpha_{k-2} \alpha_1^2 \in R^{[3]}_{BP}\bigl(2^k\bigr). $$ 
The element $\alpha_{k-2} \alpha_1^2$ detects an element of the homotopy root
invariant $R(2^k)$.

\item[${k \equiv 4 \pmod 4}$] The element $\alpha_k$ is a permanent
cycle in the ANSS, and detects the element of order $2$ in the image of $J$
in the $(2k{-}1)$--stem.  This is the root invariant $R\bigl(2^k\bigr)$.
\end{description}

\section{The first two $BP$--filtered root invariants of
$\alpha_{i/j}$}\label{sec:v1filt}
In this section we will compute
$R^{[k]}_{BP}(\alpha_{i/j})$ for $k = 1,2$ 
using the indeterminacy calculations of \fullref{sec:indeterminacy}.
We will then analyze the higher $BP$--filtered root invariants using the
theorems of \fullref{sec:results}.

The following proposition gives 
the first few filtered root invariants of
the $\alpha_{i/j}$.  Since we actually use the
homotopy root invariant to determine these filtered root invariants, this
proposition gives no new information.  
However, we do see the indeterminacy group
$A_{i/j}$ adding essential terms to the $BP\wedge \td{BP}$--root invariant.

\begin{prop}
The low dimensional filtered root invariants of the elements
$\alpha_{i/j}$ 
are given (up to a multiple in $\ZZ_{(2)}^\times$) by the following table.
\begin{center}
\begin{tabular}{|c|c|c|}
\hline
$x$ & $R_{BP}^{[1]}(x)$ & $R_{BP}^{[2]}(x)$ \\
\hline
$\alpha_1$ & $\alpha_{2/2}$ & $-$ \\
$\alpha_{2/2}$ & $\alpha_{4/4}$ & $-$ \\
$\alpha_2$ & $v_1\alpha_{4/4}$ & $\alpha_{4/4}\alpha_1$ \\
\hline
\end{tabular}
\end{center}
\end{prop}

\begin{proof}
\fullref{prop:indeterminacy}, combined with
\fullref{prop:BP^1root}, gives the following values for $R^{[1]}_{BP}$.
\begin{align*}
R^{[1]}_{BP}(\alpha_1) & \doteq \wtilde{\beta}_1 + c_1\cdot v_1\alpha_1 \\
R^{[1]}_{BP}(\alpha_{2/2}) & \doteq \wtilde{\beta}_{2/2} + c_2 \cdot x_7 \\
R^{[1]}_{BP}(\alpha_{2}) & \doteq \wtilde{\beta}_2 + c_3 \cdot v_1 x_7.
\end{align*}
\fullref{thmB} implies that 
\begin{equation}\label{eq:thmB1}
d_1\bigl(R^{[1]}_{BP}(\alpha_{i/j})\bigr) = 2\cdot R^{[2]}_{BP}(\alpha_{i/j})
\end{equation}
for the values of $i$ and $j$ we are considering.  The known root
invariants (see Mahowald--Ravenel \cite{MahowaldRavenel}) of
these elements are
$$R(\alpha_1) = R(\eta) \doteq \nu \qquad
R(\alpha_{2/2}) = R(\nu) \doteq \sigma \qquad
R(\alpha_2) = R(2\nu) \doteq \sigma\eta.$$
In the ANSS we have the following representatives of 
elements of the $E_2$--term.
\begin{align*}
\alpha_{2/2} & = \wtilde{\beta}_1 \\
\alpha_{4/4} & \equiv \wtilde{\beta}_{2/2} + x_7 \pmod 2
\end{align*}
We also have the following $d_1$--differentials in the ANSS.
\begin{align*}
d_1(x_7) &\equiv d_1(\wtilde{\beta}_{2/2}) \equiv 2\beta_{2/2} & \pmod 4 \\ 
d_1(\wtilde{\beta}_2) &\equiv 2\beta_2 & \pmod 4 \\
d_1(\wtilde{\beta}_2+v_1x_7) &\equiv 2\alpha_1\alpha_{4/4} & \pmod 4
\end{align*}
The only way these differentials can be compatible with
Equation~\eqref{eq:thmB1} and \fullref{thmA} is for the coefficients
$c_i$ to have the following values.
$$c_1 \equiv 0 \pmod 2, \qquad
c_2 \equiv 1 \pmod 2, \qquad
c_3 \equiv 1 \pmod 2.\proved$$
\end{proof}

The second filtered root invariants of the rest of the $\alpha_{i/j}$'s are
given by the following proposition.

\begin{prop}
The filtered root invariant $R_{BP}^{[2]}(\alpha_{i/j})$ for $3i-j \ge 6$ 
contains (up to a multiple in $\ZZ_{(2)}^\times$) the element
$$
\beta_{i/j} + 
\begin{cases}
c\cdot \alpha_1 \wtilde{\alpha}_{3i-j-1} & \qquad \text{ $j$ odd}  
\\
0 & \qquad \text{$j$ even}
\end{cases}
$$
with $c \in \ZZ_{(2)}$.  Here $\wtilde{\alpha}_{k}$ represents the ANSS element
$\alpha_{k/l}$ with $l$ maximal.
\end{prop}

\begin{proof}
\fullref{prop:indeterminacy}, together with
\fullref{prop:BP^1root} gives
$$ \wtilde{\beta}_{i/j} + c \cdot x \dotin R^{[1]}_{BP}(\alpha_{i/j}) $$
where $c$ is an element of $\ZZ_{(2)}$ and $x \in BP_*\td{BP}$ has the 
property that the Adams--Novikov differential $d_1(x)$ is given by
\begin{equation}\label{eq:d1x}
d_1(x) \doteq
\begin{cases}
2 \alpha_1 \wtilde{\alpha}_{3i-j-1}, & \qquad \nu_2(i) \le 3i-j-2, \text{ $j$ odd}
\\
0 & \qquad \text{otherwise.}
\end{cases}
\end{equation}
We claim that the condition $\nu_2(i) \le 3i - j - 2$ is always satisfied for
$3i - j \ge 6$ where $i$ and $j$ are such that $\alpha_{i/j}$ exists in the
Adams--Novikov $E_2$--term.  Indeed, for $\alpha_{i/j}$ to exist we must have 
$j \le \nu_2(i) + 2$, from which it follows that 
$$ 3i - \nu_2(i) - 4 \le 3i - j - 2. $$
Therefore it suffices to show that $\nu_2(i) \le 3i - \nu_2(i) - 4$, or
equivalently $2\nu_2(i) \le 3i - 4$.  The latter is true for $i \ge 2$, and
the condition $3i-j \ge 6$ in particular implies that $i \ge 2$.
Thus the condition in the first case of Equation \eqref{eq:d1x} may be
simplified to simply read ``$j$ odd.''

There are Adams--Novikov differentials 
$$ d_1(\wtilde{\beta}_{i/j}) \equiv 2\beta_{i/j} \pmod 4 . $$
\fullref{thmB} applies to give
$$ d_1\bigl(R^{[1]}_{BP}(\alpha_{i/j})\bigr) =
2\cdot R^{[2]}_{BP}(\alpha_{i/j}) $$
and the result follows.
\end{proof}

\section{Higher $BP$ and $H\FF_2$--filtered root invariants of
some $v_1$--periodic elements}\label{sec:froots}

We denote the Eilenberg--MacLane spectrum $H\FF_2$ by
$H$.  We describe what happens first in the ASS, and then the ANSS.  All
of the algebraic root invariants used were taken from Bruner \cite{Bruner}.  These
algebraic root invariants coincide with the first non-trivial filtered 
root invariants by 
\fullref{thmD}.  The
computations are summarized below.  We will show in
\fullref{sec:hroots} that in each of these cases the top filtered root
invariant successfully detects the homotopy root invariant through the use
of \fullref{thmA}.

\begin{tabbing}
$R(\wbar{\alpha}_1)$ \hspace{.3in} \= 
\parbox[t]{4.40in}{ASS:
The element $\wbar{\alpha}_1 = \eta$ is detected by $h_1$.  
We have $h_2 \in R_{\alg}(h_1)$.  The element $h_2$ detects $\nu$.}
\\
\> \parbox[t]{4.40in}{ANSS:
We have $\alpha_{2/2} \in R^{[1]}_{BP}(\wbar{\alpha}_1)$.  The element
$\alpha_{2/2}$ also detects $\nu$.}
\\
$R(\wbar{\alpha}_1^2)$ \>
\parbox[t]{4.40in}{ASS: 
The algebraic root invariant is given 
by $h_2^2 \in R_{\alg}(h_1^2)$,
and $h_2^2$ detects $\nu^2$.}
\\
$R(\wbar{\alpha}_1^3)$ \>
\parbox[t]{4.40in}{ASS:
We have $h_2^3 \in R_{\alg}(h_1^3)$, which detects $\nu^3$.}
\\
$R(\wbar{\alpha}_{2/2})$ \> 
\parbox[t]{4.40in}{ASS:
The element $\wbar{\alpha}_{2/2}$ is detected by $h_2$.  The algebraic root
invariant is given by $h_3 \in R_{\alg}(h_2)$. The element $h_3$ detects 
$\sigma$.}  
\\
\> \parbox[t]{4.40in}{ANSS:
The first filtered root
invariant is given by $\alpha_{4/4} \in R^{[1]}_{BP}(\alpha_{2/2})$, and
$\alpha_{4/4}$ detects $\sigma$.}
\\
$R(\wbar{\alpha}_2)$ \> 
\parbox[t]{4.40in}{ASS:
The element $\wbar{\alpha}_2$ is detected by $h_0 h_2$, and $h_1
h_3 \in R_{\alg}(h_0 h_2)$ detects $\eta \sigma$.}
\\
\> \parbox[t]{4.40in}{ANSS:
The second filtered root invariant is given by $\alpha_{4/4}\alpha_1 
\in R^{[2]}_{BP}(\wbar{\alpha}_2)$, which detects $\eta\sigma$.}
\\
$R(\wbar{\alpha}_{4/4})$ \> 
\parbox[t]{4.40in}{ASS: 
The element $\wbar{\alpha}_{4/4}$ is detected by $h_3$ in the ASS.  The
algebraic root invariant is given by 
$h_4 \in R_{\alg}(h_3)$, which coincides with the first
filtered root invariant $R^{[1]}_{H}(h_3)$ by \fullref{thmD}.
The Hopf invariant $1$
differential $d_2(h_4) = h_0 h_3^2$ allows one apply \fullref{thmB} and
get $h_3^2 \in R^{[2]}_{H}(h_3)$.  The element $h_3^2$ detects $\sigma^2$.}
\\
\> \parbox[t]{4.40in}{ANSS:
We have the second filtered root invariant 
$\beta_{4/4} \in R^{[2]}_{BP}(\wbar{\alpha}_{4/4})$, and $\beta_{4/4}$ detects
$\sigma^2$.}
\\
$R(\wbar{\alpha}_{4/3})$ \>
\parbox[t]{4.40in}{ASS:
The element $\wbar{\alpha}_{4/3}$ is detected by
$h_0 h_3$ in the ASS.  The algebraic root invariant is $h_1 h_4 \in 
R_{\alg}(h_0 h_3)$, which detects
$\eta_{4}$.}
\\
\> \parbox[t]{4.40in}{ANSS:
The second filtered root invariant is
$\beta_{4/3} + c \alpha_{8/5} \alpha_1 \in R^{[2]}_{BP}(\wbar{\alpha}_{4/3})$, 
for
$c \in \ZZ/2$.
The element $\beta_{4/3}$ detects $\eta_4$.  The value of $c$ is irrelevant
--- the AHSS differentials of Mahowald \cite{Mahowald} imply that if $\eta_4$ is in
the root invariant $R(2\sigma)$ (which it is), 
then $\alpha_1\alpha_{8/5}$ is in the
indeterminacy of this root invariant.  In the AHSS for
$\pi_*(P_{-10}^\infty)$ the element $\wbar{\alpha}_1\wbar{\alpha}_{8/5}[-10]$ 
is the target of a differential supported by $\wbar{\alpha}_5[-2]$.}
\\
$R(\wbar{\alpha}_{4/2})$ \> 
\parbox[t]{4.40in}{ASS:
The element $\wbar{\alpha}_{4/2}$ is detected by 
$h_0^2 h_3$.
The algebraic root invariant is 
$h_1^2 h_4 \in R_{\alg}(h_0^2 h_3)$, and this detects $\eta \eta_4$.} 
\\
\> \parbox[t]{4.40in}{ANSS:
The second filtered root invariant is given by
$\beta_{4/2} \in R^{[2]}_{BP}(\wbar{\alpha}_{4/2})$,
where $c \in \ZZ/2$.
The element 
$\beta_{4/2}$
is a permanent cycle but does not survive to the homotopy root invariant.
We appeal to \fullref{thmC} to find a higher filtered root invariant. 
In $P^{-12}_{-13}$ the $-12$--cell attaches 
to the
$-13$--cell with degree $2$ attaching map, and the $-11$--cell attaches to
the $-13$--cell in $P^{-11}_{-13}$ by $\eta$.  There is a hidden extension in 
the ANSS given by
$\wbar{\alpha}_1 \cdot \wbar{\alpha}_1 \wbar{\beta}_{4/3} = 
2 \cdot \wbar{\beta}_{4/2}$, which indicates that
we have a higher filtered root invariant $\alpha_1 \beta_{4/3} \in
R^{[3]}_{BP}(\wbar{\alpha}_{4/2})$.  The element $\alpha_1 \beta_{4/3}$ detects 
$\eta \eta_4$.}
\\
$R(\wbar{\alpha}_{4})$ \> 
\parbox[t]{4.40in}{ASS:
The element $\wbar{\alpha}_{4}$ is detected by $h_0^3 h_3$, with algebraic
root invariant $h_1^3 h_4
\in R_{\alg}(h_0^3 h_3)$ which detects $\eta^2 \eta_4$.}
\\
\> \parbox[t]{4.40in}{ANSS:
The second filtered root invariant is given by 
$(\beta_4 + c\alpha_1\alpha_{10} \in R^{[2]}_{BP}(\wbar{\alpha}_4)$.  
It turns out that $\alpha_1\alpha_{10}$ is in the indeterminacy of 
$R^{[2]}_{BP}(\wbar{\alpha}_4)$, since in the AHSS for $P_{-14}^\infty$, we
have $d_2(\alpha_{10}[-12]) = \alpha_{10}\alpha_1[-14]$.  Therefore we may
as well set $c = 0$.
The element $\beta_4$ corresponds 
to the element $g$
in the ASS.  
The element $\beta_4$ is a permanent
cycle, so we use \fullref{thmC} to look for a higher root invariant.  

In $P^{-12}_{-15}$ the $-14$--cell attaches to the $-15$--cell by the degree
$2$ map, and there is a Toda bracket
$\bigl\langle{P^{-12}_{-15}}\bigr\rangle(-) \capeq
\bra{2, \alpha_1, -}$.

There is a hidden extension in the ANSS given by
$2 \wbar{\beta}_4 = (\wbar{\beta}_{4/4}\wbar{\beta}_{2/2})/2$, and in the ANSS
there is a Toda bracket $(\beta_{4/4}\beta_{2/2})/2 \in \bra{2, \alpha_1,
\beta_{4/3} \alpha_1^2}$.
}
\\
\> \parbox[t]{4.40in}{

We conclude using \fullref{thmC} that we
have the higher filtered root invariant
$\beta_{4/3} \alpha_1^2 \in R^{[4]}_{BP}(\wbar{\alpha}_4)$.  The element 
$\beta_{4/3} \alpha_1^2$ detects $\eta^2 \eta_4$.}
\\
$R(\wbar{\alpha}_{4/4} \wbar{\alpha}_1)$ \>
\parbox[t]{4.40in}{ASS:
The element $\wbar{\alpha}_{4/4} \wbar{\alpha}_1$ is detected by $h_3h_1$
with algebraic root invariant $h_4h_2 \in R_{\alg}(h_3h_1)$ which detects
$\nu^*$.}
\\
\> \parbox[t]{4.40in}{ANSS:
The root invariant of $\wbar{\alpha}_{4/4}$ will turn out to be given by 
$\wbar{\beta}_{4/4}$, 
so one
might initially suspect that $\wbar{\alpha}_{2/2} \wbar{\beta}_{4/4}$ 
would detect the
homotopy root invariant $R(\wbar{\alpha}_{4/4} \wbar{\alpha}_1)$.  However, the
$-6$--cell attaches to the $-10$--cell with attaching map $\nu$.
Thus in the AAHSS for $P_{-10}$, there is a differential
$d_4(\beta_{4/4}[-6]) = \alpha_{2/2} \cdot \beta_{4/4}[-10]$.  Looking at
the attaching map structure of $P^{-6}_{-11}$, this differential actually
pushes the root invariant to $\beta_{4/2,2} \in 
\bra{\alpha_{2/2}, 2, \beta_{4/4}}$, and this element detects $\nu^*$.}
\\
$R(\wbar{\alpha}_{4/4} \wbar{\alpha}_1^2)$ \>
\parbox[t]{4.40in}{ASS:
The element $\wbar{\alpha}_{4/4} \wbar{\alpha}_1^2$ is detected by $h_3 h_1^2$,
with algebraic root invariant $h_4 h_2^2 \in R_{\alg}(h_3 h_1^2)$.  
which detects $\nu^* \nu$. }
\\
$R(\wbar{\alpha}_5)$ \> 
\parbox[t]{4.40in}{ASS:
The element $\wbar{\alpha}_5$ is detected by $Ph_1$, with algebraic root
invariant $h_2 g
\in R_{\alg}(Ph_1)$ which detects
$\nu \wbar{\kappa}$. } 
\\
\> \parbox[t]{4.40in}{ANSS:
The second filtered root invariant is given by
$\beta_5 + c\alpha_{13} \alpha_1 \in R^{[2]}_{BP}(\wbar{\alpha}_5)$, 
where $c \in
\ZZ/2$.  \fullref{thmB} applies to the differential
$d_3(\beta_5 + c \alpha_{13}\alpha_1) = 
\alpha_1^2 \eta_{3/2}$ 
to give the higher filtered root invariant $\alpha_1 \eta_{3/2} \in 
R^{[4]}_{BP}(\wbar{\alpha}_5)$.
The element $\alpha_1 \eta_{3/2}$ corresponds to the element $h_4c_0h_1$ in
the ASS.  Using the May spectral sequence, we see that there is a Massey
product $h_2^2g \in \bra{h_0, h_1, h_4c_0h_1}$.  The element $h_2^2g$ is a
permanent cycle which corresponds to the element $\Pi \beta_{2/2}$ in the
ANSS.  We conclude that in the ANSS there is a hidden Toda bracket
$\wbar{\alpha}_{2/2} \cdot \wbar{\beta}_4 \wbar{\alpha}_1^3/8 = \br{\Pi
\beta_{2/2}} \in \bra{2,
\wbar{\alpha}_1, \wbar{\alpha}_1 \wbar{\eta}_{3/2}}$.  In $P^{-15}_{-19}$, 
the $-16$--cell attaches to the $-19$--cell with the Toda bracket $\bra{2,
\wbar{\alpha}_1, -}$ and the $-15$--cell attaches to the $-19$--cell with
attaching map $\wbar{\alpha}_{2/2}$.  We conclude, using \fullref{thmC},
that we have another higher filtered root invariant $\beta_4 \alpha_1^3/8
\in R^{[5]}_{BP}(\wbar{\alpha}_5)$.  The element $\beta_4 \alpha_1^3/8$ detects
$\nu \wbar{\kappa}$.}
\\
$R(\wbar{\alpha}_5\wbar{\alpha}_1)$ \>
\parbox[t]{4.40in}{ASS:
The element $\wbar{\alpha}_5\wbar{\alpha}_1$ is detected by $Ph_1^2$, with
algebraic root invariant $r \in R_{\alg}(Ph_1^2)$.  \fullref{thmD} implies
that we have the filtered root invariant $r \in R^{[6]}_{H}(
\wbar{\alpha}_5\wbar{\alpha}_1)$.
The element $r$ supports
the Adams differential $d_3(r) = h_1 d_0^2$.
Thus, we may use \fullref{thmB}
to deduce that there is a higher filtered root invariant
$d_0^2 \in R^{[8]}_{H}(
\wbar{\alpha}_5\wbar{\alpha}_1)$ which detects $\kappa^2 = \epsilon\wbar{\kappa}$.}
\\
$R(\wbar{\alpha}_5\wbar{\alpha}_1^2)$ \>
\parbox[t]{4.40in}{ASS: 
The element $\wbar{\alpha}_5\wbar{\alpha}_1^2$ is detected by $Ph_1^3$,
with algebraic root invariant $h_1q \in R_{\alg}(Ph_1^3)$.  }
\\
$R(\wbar{\alpha}_{6/2})$ \> 
\parbox[t]{4.40in}{ASS:
This element presents a very interesting story: it is an instance
where the ASS seems to give us nothing yet the ANSS tells us the homotopy root
invariant.  The element $\wbar{\alpha}_{6/2}$ corresponds to the element 
$Ph_2$ in the ASS,
with algebraic root invariant $h_0^3 h_4^2 \in R_{\alg}(Ph_2)$.  The element 
$h_0^3 h_4^2$ is killed by a
differential in the ASS, and it is not clear what a candidate for the root
invariant should be.}
\\
\> \parbox[t]{4.40in}{ANSS:
The second filtered root invariant is given by
$\beta_{6/2} \in R^{[2]}_{BP}(\wbar{\alpha}_{6/2})$.  There is an Adams--Novikov
differential $d_5 (\beta_{6/2}) = \alpha_{2/2}
\cdot \Pi \beta_{2/2}$.  \fullref{thmB} indicates that we have the
higher filtered root invariant $\Pi \beta_{2/2}
\in R^{[6]}_{BP}(\alpha_{6/2})$ which detects $\nu^2 \wbar{\kappa}$.}
\\
$R(\wbar{\alpha}_{6})$ \>
\parbox[t]{4.40in}{ASS:
The element $\wbar{\alpha}_{6}$ is detected by $Ph_2h_0$, with algebraic root
invariant $q \in R_{\alg}(Ph_2h_0)$.}
\end{tabbing}

\section{Homotopy root invariants of some $v_1$--periodic 
elements}\label{sec:hroots}

In this section we will use the filtered root invariant computations of
\fullref{sec:froots} to compute some homotopy root invariants.  These
computations are summarized in the following theorem.  The first few are
well-known (see Mahowald--Ravenel \cite{MahowaldRavenel}).

\begin{thm}\label{thm:hroots}
We have the following table of homotopy root invariants (up to some
multiple in $\ZZ_{(2)}^\times$).

\begin{center}
\begin{tabular}{|c|c|}
\hline
$x$ & $R(x)$ \\
\hline
$\wbar{\alpha}_1 = \eta$ & $\nu$ \\
$\wbar{\alpha}_1^2 = \eta^2$ & $\nu^2$ \\
$\wbar{\alpha}_1^3 = \eta^3$ & $\nu^3$ \\
$\wbar{\alpha}_{2/2} = \nu$ & $\sigma$ \\
$\wbar{\alpha}_{2} = 2\nu$ & $\sigma\eta$ \\
$\wbar{\alpha}_{4/4} = \sigma$ & $\sigma^2$ \\
$\wbar{\alpha}_{4/3} = 2\sigma$ & $\eta_4$ \\
$\wbar{\alpha}_{4/2} = 4\sigma$ & $\eta\eta_4$ \\
\hline
\end{tabular}
\qquad
\begin{tabular}{|c|c|}
\hline
$x$ & $R(x)$ \\
\hline
$\wbar{\alpha}_{4} = 8\sigma$ & $\eta^2\eta_4$ \\
$\wbar{\alpha}_{4/4}\wbar{\alpha}_1 = \eta\sigma$ & $\nu^*$ \\
$\wbar{\alpha}_{4/4}\wbar{\alpha}_1^2 = \sigma\eta^2$ & $\nu\nu^*$ \\
$\wbar{\alpha}_5$ & $\nu\wbar{\kappa}$ \\
$\wbar{\alpha}_5\wbar{\alpha}_1$ & $\kappa^2 = \epsilon\wbar{\kappa}$ \\
$\wbar{\alpha}_5\wbar{\alpha}_1^2$ & $\eta\wbar{q}$ \\
$\wbar{\alpha}_{6/2}$ & $\nu^2\wbar{\kappa}$ \\
$\wbar{\alpha}_6$ & $\wbar{q}$ \\
\hline
\end{tabular}
\end{center}

Some of these root invariants may have indeterminacy.
\end{thm}

We pause to remark on the elements that show up as root invariants in this
table.  With the exception of $\nu$, $\sigma$, $\sigma\eta$, and $\nu^3 =
\sigma\eta^2$,
all of these elements are
$v_2$--periodic.  This was shown by Mahowald in \cite{Mahowaldv2},
and Hopkins and Mahowald in \cite{HopkinsMahowald}.  We remind the reader
that due to an error in Davis--Mahowald \cite{DavisMahowald}, one must replace $8k$ with
$32k$ in \cite{Mahowaldv2}.  Recently, Hopkins and Mahowald have produced
some $v_2^{32}$--self maps \cite{HopkinsMahowaldself}.  In particular, in
\cite[Problem~4]{Mahowald} (see also \cite{DavisMahowald}), a list of
$v_2$--periodic elements in $\pi_*(S)$ are given which are the first few
homotopy Greek letter elements $\beta_i^h$:
$$ \nu, \nu^2, \nu^3, \eta^2\eta_4, \nu\wbar{\kappa}, \epsilon\wbar{\kappa},
\eta\wbar{q}, \ldots $$
Our computations show that these elements appear as the iterated root
invariants $R(R(2^i))$ for $i \le 7$.

The rest of this section is devoted to proving \fullref{thm:hroots}.
We use \fullref{corA} to deduce the homotopy root invariants from our
filtered root invariants.  In our range, the first
part of
\fullref{corA} is easier to check using the ANSS rather than the ASS, 
since
there are fewer elements to check.  However, we have only determined the
$BP$--filtered root invariants for the elements in
Adams--Novikov filtration $1$.  For the $v_1$--periodic elements in higher
Adams--Novikov filtration, we must use our $H\FF_2$--filtered root 
invariants and the
ASS.  Since the author's knowledge of the $2$--primary ANSS does not include
$\beta_6$, we also use the ASS to compute $R(\alpha_6)$.
The second part of \fullref{corA} is verified afterwards.

Our computations for the first part of \fullref{corA} using 
$BP$--filtered root invariants and the ANSS
are summarized in the following table.

\begin{center}
\begin{tabular}{|c|c|c|c|l|}
\hline
$x$ & $R^{[k]}_{BP}(x)$ & $-N$ & $\{\gamma_i[n_i]\}$ & diff \\
\hline
$\wbar{\alpha}_1$ & $\alpha_{2/2}$ & $-3$ &
 &  \\
$\wbar{\alpha}_{2/2}$ & $\alpha_{4/4}$ & $-5$ & 
$\beta_{2/2}[-4]$ & $\leftarrow \alpha_1[2]$ \\
$\wbar{\alpha}_2$ & $\alpha_{4/4}\alpha_1$ & $-6$ &
&  \\
$\wbar{\alpha}_{4/4}$ & $\beta_{4/4}$ & $-8$ &
 &  \\
$\wbar{\alpha}_{4/3}$ & $\beta_{4/3}$ & $-10$ &
$\beta_{4/4}\alpha_1[-9]$ & $\leftarrow \beta_{4/4}[-7]$ \\
$\wbar{\alpha}_{4/2}$ & $\beta_{4/3}\alpha_1$ & $-11$ &
 &  \\
$\wbar{\alpha}_{4}$ & $\beta_{4/3}\alpha_1^2$ & $-12$ &
 &  \\
$\wbar{\alpha}_5$ & $\alpha_{1}^3\beta_4/8$ & $-15$ & 
 &  \\
$\wbar{\alpha}_{6/2}$ & $\Pi \beta_{2/2}$ & $-16$ &
 &  \\
\hline
\end{tabular}
\end{center}

We explain the columns of the table:
\begin{center}
\begin{tabular}{ll}
$x$ & The element we wish to compute the root invariant of. \\
$R^{[k]}_{BP}(x)$ & The top $BP$--filtered root invariant of $x$, 
which we want to show \\
& detects $R(x)$ using \fullref{corA}. \\
$-N$ & The cell of $P_{-\infty}^\infty$ which carries the 
filtered root invariant.
\end{tabular}
\end{center}

\begin{center}
\begin{tabular}{ll}
$\{\gamma_i\lbrack n_i \rbrack\}$ & The collection of 
elements in the $E_1$--term of
the AAHSS \\
& converging to $\ext(BP_*P_{-N})$  
which could detect the difference \\
& between the top filtered root invariant and the homotopy root invariant. \\
& We exclude elements of infinite $\alpha_1$--towers, since these are
  always the \\
& source or target of a $d_2$--differential in the AAHSS. \\
diff & Each element $\gamma_i[n_i]$ turns out to be ineligible to
detect the difference, \\
& since it is the target the indicated AAHSS differential.
\end{tabular}
\end{center}

We now deal with the leftover elements using the $H\FF_2$--filtered root
invariants and the ASS.  We have the following 
$H\FF_2$--filtered root invariants
which are permanent cycles in the ASS.

\begin{center}
\begin{tabular}{|c|c|c|}
\hline
$x$ & $R^{[k]}_{H}(x)$ & $-N$ \\
\hline
$\wbar{\alpha}_{4/4}\wbar{\alpha}_1$ & $h_4h_2$ & $-11$  \\
$\wbar{\alpha}_{4/4}\wbar{\alpha}_1^2$ & $h_4h_2^2$ & $-13$ \\
$\wbar{\alpha}_5\wbar{\alpha}_1$ & $d_0^2$ & $-19$ \\
$\wbar{\alpha}_5\wbar{\alpha}_1^2$ & $h_1q$ & $-23$ \\
$\wbar{\alpha}_6$ & $q$ & $-22$ \\
\hline
\end{tabular}
\end{center}

Using \fullref{corA}, we see that to verify that these filtered root
invariants detect
homotopy root invariants, we must first check that there are no elements of
$\pi_{t-1}(P_{-N})$ of Adams filtration greater than $k$ which can detect
the root
invariant on a higher cell.  
We handle this on a case-by-case basis, with the aid of the
computations of Mahowald in \cite{Mahowaldmetastable}, and the computer
$\ext$ computations of Bruner \cite{Brunerweb}.
In the following analysis,
we omit the elements detected in the AHSS by 
$v_1$--periodic elements.  
These elements cannot be root invariants in the stems we are considering 
by the following lemma.

\begin{lem}\label{lem:Jgamma_i}
Suppose that $\gamma[n] \in \pi_j(S^n)$ is an element of the $E_1$--term of
the AHSS for $\pi_j(P_{-\infty}^\infty)$ where $\gamma$ is a $v_1$--periodic 
element.  Then $\gamma[n]$ is either the source or target of a non-trivial
AHSS differential unless we are in one of the following cases (in which
case we do not know whether $\gamma[n]$ is the source or target of a
differential).
\begin{itemize}
\item $\gamma = \nu$ and $j \equiv 0 \pmod 4$
\item $\gamma = B_i$ with $j \equiv 6 \pmod 8$ and $i \le \nu_2(j+2)-2$
\item $\gamma = \sigma$, $\sigma\eta$, or  $\sigma\eta^2$ with $j \equiv 6
\pmod 8$ (The behavior of these elements is slightly anomalous, due to the
presence of $\nu^2$, $\epsilon$, and $\eta\epsilon$.)
\end{itemize}
Here $B_i$ is $i^\mathrm{th}$ generator of the image of the
$J$--homomorphism. Therefore, if $i = 4a+b$ with $0 \le b \le 3$, we have
$$ B_i = \begin{cases}
\wtilde{\alpha}_{4a}\wbar{\alpha}_1^2, & b = 0 \\
\wtilde{\alpha}_{4a+2}, & b = 1 \\
\wtilde{\alpha}_{4(a+1)}, & b = 2 \\
\wtilde{\alpha}_{4(a+1)}\wbar{\alpha}_1, & b = 3
\end{cases}
$$
where $\wtilde{\alpha}_k$ is the element $\wbar{\alpha}_{k/l}$ with $l$ maximal.
\end{lem}

\begin{proof}
Mahowald \cite[Theorem~4.6]{MahowaldDescriptive} states that you can lift the
differentials from the $J$--homology modified AHSS to the (double suspension)
EHP spectral
sequence.  The proofs of the announcements in \cite{MahowaldDescriptive}
are the subject of \cite{Mahowald}.  Since
the EHP spectral sequence maps to the AHSS for $\pi_*(P^\infty)$, the
differentials in the AHSS follow.  We then get the result for the AHSS for
$P^\infty_{-\infty}$ by transporting our differentials with James periodicity.
\end{proof}

The first part of \fullref{corA} on the remaining elements is verified 
as follows. 

\noindent
$\wbar{\alpha}_{4/4} \wbar{\alpha}_1$:  
According to the tables of Mahowald \cite{Mahowaldmetastable}, the only elements of
$\pi_{7}(P_{-11})$ of Adams filtration greater than $2$ are $v_1$--periodic.

\noindent
$\wbar{\alpha}_{4/4} \wbar{\alpha}_1^2$:
According to the tables of \cite{Mahowaldmetastable}, there are no elements of
$\pi_{8}(P_{-13})$ of Adams filtration greater than $3$.

\noindent
$\wbar{\alpha}_{5} \wbar{\alpha}_1$:
Examining the tables of \cite{Mahowaldmetastable}, the only elements of
$\pi_{9}(P_{-19})$ of Adams filtration greater than $8$ have trivial image
in $\pi_9(P_{-18})$.  Therefore, none of them can detect a root invariant
carried by a cell above the $-19$ cell.

\noindent
$\wbar{\alpha}_{5} \wbar{\alpha}_1^2$:
Examining the tables of Bruner \cite{Brunerweb}, we find the following pattern of
generators in $\ext(H_*P_{-23})$.
\begin{center}
\labellist\tiny
\pinlabel {$s$} [r] at 13 181
\pinlabel {$7$} [r] at 6 63
\pinlabel {$10$} [t] at 72 11
\pinlabel {$t{-}s$} [tl] at 117 18
\pinlabel {$P^2d_0h_1^2[-23]$} [b] at 47 168
\pinlabel {$P^2e_0[-23]$} [l] at 76 137
\pinlabel {$h_1q[-23]$} [t] at 37 50
\pinlabel {\small $\ext(P_{-23})$} [tl] at 105 196
\endlabellist
\includegraphics[width=1.5in]{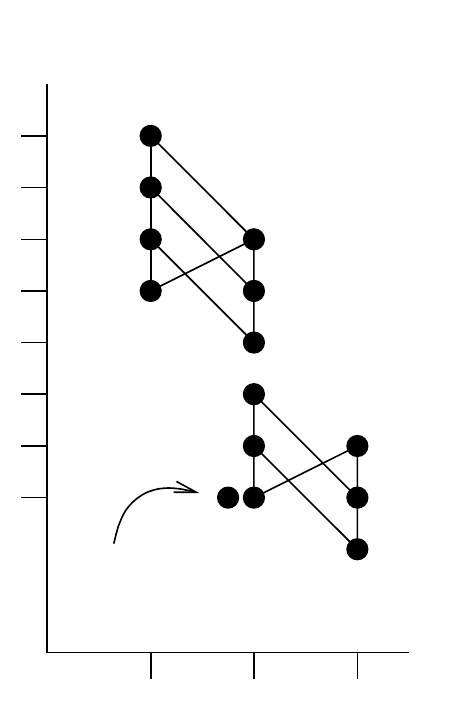} 
\end{center}
Some of the elements are labeled with their AAHSS names.  These names were
deduced from the AAHSS differentials computed in \cite{Mahowaldmetastable}.
The Adams differentials originating from the elements in Adams filtration
$6$ and $7$ may be deduced by extrapolating differentials computed in
\cite{Mahowaldmetastable} using $h_0$, $h_1$, and $h_2$ multiplication.  The
inclusion of the bottom cell
$$ S^{-23} \rightarrow P_{-23} $$
induces the differential $d_2(P^2e_0[-23]) = P^2d_0h_1^2[-23]$.  The 
rest of the differentials are then forced by $h_0$ multiplication.  We
deduce that the only elements of $\pi_{10}(P_{-23})$ of Adams filtration
greater than $7$ are the $v_1$--periodic elements.

\noindent
$\wbar{\alpha}_{6}$:
From the tables of Bruner \cite{Brunerweb}, we find the following portion of
$\ext(H_*P_{-22})$.
\begin{center}
\labellist\tiny
\pinlabel {$s$} [r] at 14 211
\pinlabel {$6$} [r] at 8 92
\pinlabel {$10$} [t] at 73 11
\pinlabel {$t{-}s$} [t] at 134 18
\pinlabel {$P^2h_2^3[-16]$} [t] at 41 165
\pinlabel {$k[-19]$} [r] at 74 104
\pinlabel {$q[-22]$} [tr] at 74 93
\pinlabel {$p[-22]$} [r] at 103 62
\pinlabel {$h_4^2h_1[-20]$} [t] at 103 47
\pinlabel {$d_1[-21]$} [l] at 113 63
\pinlabel {$h_0l[-15]$} [bl] at 81 122
\pinlabel {\small $\ext(P_{-22})$} [t] at 122 215
\endlabellist
\includegraphics[width=1.5in]{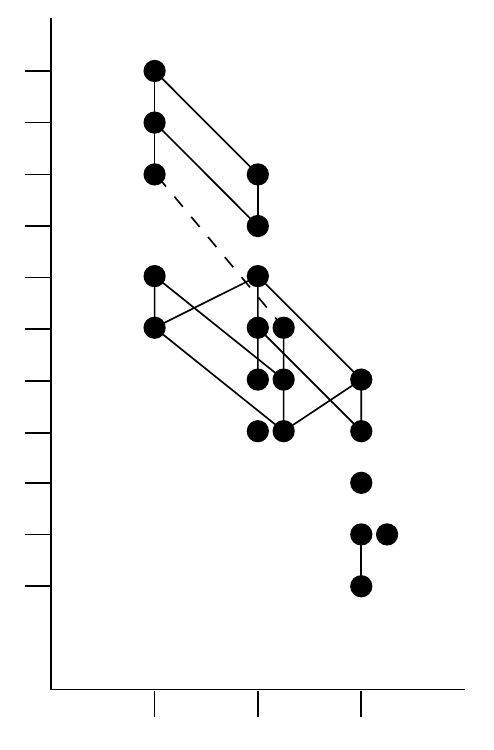}
\end{center}
All of the $d_2$ differentials shown are extrapolated from differentials in
the charts of Mahowald \cite{Mahowaldmetastable} using $h_0$, $h_1$, and
$h_2$--multiplication.  The remaining two elements that could detect
elements of higher Adams filtration, $h_0l[-15]$ and $k[-19]$, must be
handled with care.  We make the following claims, which combine to show
there are no classes in $\pi_{10}(P_{-22})$ of Adams filtration greater
than $6$ which could detect the root invariant of $\wbar{\alpha}_6$.
\begin{enumerate}
\item There is an Adams differential $d_3(h_0l[-22]) = P^2h_2^3[-16]$ (as
indicated by a dashed line in the chart).
\item The element $k[-19]$ is a non-trivial permanent cycle which detects
an element $\gamma \in \pi_{10}(P_{-22})$.
\item The image of $\gamma$ in $\pi_{10}(P_{-21})$ cannot agree with 
the image of $\wbar{\alpha}_6$ under the
composite $S^{-1} \rightarrow P_{-\infty}^\infty \rightarrow P_{-21}$.
\end{enumerate}

\noindent
\textbf{Proof of (1)}\qua In the ASS for $\pi_*(S^0)$, there is a differential
$$ d_2(h_0l) = Pe_0h_2^2. $$
In the AAHSS for $\ext(H_*P_{-22})$, there is a differential
$$ d_6(Pe_0h_1[-16]) = \langle Pe_0h_1 , h_2 , h_1 
\rangle [-22] =
Pe_0h_2^2[-22]. $$
However, in the ASS for $\pi_*(S^0)$, there is a differential
$$ d_2(Pe_0h_1) = P^2h_2^3. $$
We conclude that in the $E_3$--term of the ASS for $\pi_*(P_{-22})$, the
elements $Pe_0h_2^2[-22]$ and $P^2h_2^3[-16]$ have been equated.  Thus, the
element $h_0l[-22]$ must kill the element $P^2h_2^3[-16]$.

\textbf{Proof of (2)}\qua  The generator of $\ext(H_*P_{-22})$ 
in $(t-s,s) = (11,5)$ cannot support a $d_2$ killing $k[-19]$ because it
does not support non-trivial $h_0$ multiplication.  The elements $p$,
$d_1$, and
$h_4^2h_1$ of $\ext(\FF_2)$ detect 
homotopy elements $\wbar{p}$, $\wbar{d}_1$, and $\eta\theta_4$ in $\pi_*(S)$.
These elements are easily seen to extend to elements 
$\wbar{p}[-22]$, $\wbar{d}_1[-21]$, and $\eta\theta_4[-20]$
of
$\pi_{11}\bigl(P_{-22}^{-20}\bigr)$.
The elements $p[-22]$, $d_1[-21]$, and $h_4^2h_1[-20]$   
detect the images of $\wbar{p}[-22]$, $\wbar{d}_1[-21]$, and $\eta\theta_4[-20]$
under the inclusion
$$ P_{-22}^{-20} \rightarrow P_{-22}. $$
Therefore, the elements $p[-22]$, $d_1[-21]$ and $h_4^2h_1[-20]$ 
must be permanent
cycles.  There are no other elements of $\ext(H_*P_{-22})$ which can kill
$k[-19]$.

\textbf{Proof of (3)}\qua 
Let $\gamma \in \pi_{10}(P_{-22})$ be detected by the permanent cycle
$k[-19]$.
Let $\nu_N$ denote the composite
$$ \pi_*(S^{-1}) \rightarrow \pi_*(P_{-\infty}) 
\rightarrow \pi_*(P_{N}). $$
Let $\gamma'$ be the image of $\gamma$ in $\pi_{10}(P_{-21})$.
Since the element 
$k[-19]$ is non-trivial in $\ext(H_*P_{-21})$, the element $\gamma'$
has Adams filtration $7$.
We must show that $\gamma'$
cannot equal $\nu_{-21}(\wbar{\alpha}_6)$.
Let $\wbar{q} \in \pi_{32}(S^0)$ be the element detected by $q$.
Examining the ASS for $S^0$ (see Mahowald--Tangora \cite{MahowaldTangora},
and Ravenel \cite{Ravenel}), we see
that the element $\wbar{q}$ extends to an element $\delta = \wbar{q}[-22]$ in
$\pi_{10}(P_{-25})$.  
Since the algebraic root invariant of $Ph_2h_0$ is $q$, the element
$\nu_{-25}(\wbar{\alpha}_6)$ is detected in $\ext(H_*P_{-22})$ by $q[-22]$.
Since $q[-22]$ also detects $\delta$, the sum $\zeta =
\nu_{-25}(\wbar{\alpha}_6) + \delta$ has Adams filtration greater than $6$. 
Consulting the charts in \cite{Brunerweb}, we find that
there is one generator in $\ext^{s,t}(H_*P_{-25})$ for $(t-s,s) = (10, 7)$,
which is detected in the AAHSS by the element $h_1q[-23]$.  
The image of $\zeta$ in $\pi_{10}(P_{-21})$ is $\nu_{-21}(\wbar{\alpha}_6)$.
Since the image
of $h_1q[-23]$ in $\ext(H_*P_{-21})$ is zero, we deduce that
$\nu_{-21}(\wbar{\alpha}_6)$ is of Adams filtration greater than $7$.  
Thus $\nu_{-21}(\wbar{\alpha}_6)$ cannot equal $\gamma'$, 
since
$\gamma'$ has Adams filtration $6$.

We have verified the first part of \fullref{corA}.
We now must fulfill the second part of \fullref{corA}. 
Suppose that we are given an 
element $x \in \pi_*(S)$ and we want to see that $y$ is an element of 
$R(x)$ for some $y
\in \pi_*(S)$, where the root invariant is carried by the $-N$--cell.
Suppose that $y$ was detected by a filtered root invariant
and that we have
already verified the first part of \fullref{corA}.
Then the root invariant of $x$ is $y$ if we can show
that the image of the element $y$ 
under the inclusion of the bottom cell
$$ \pi_*(S^{-N}) \rightarrow \pi_*(P_{-N}) $$
is nontrivial.  

\begin{lem}
Let $E$ be either $H\FF_2$ or $BP$ and suppose that $\wtilde{y} \in \ext(E_*)$
detects $y$ in the $E$--ASS.
It suffices to show the following two things:
\begin{enumerate}
\item The element $\wtilde{y}[-N]$ is not the target of a differential in the
AAHSS.

\item This element of $\ext(E_*P_{-N})$ which $\wtilde{y}[-N]$ detects 
is not the target of an $E$--ASS differential.
\end{enumerate}
\end{lem}

We have the following convenient proposition.

\begin{prop}
If $z$ is a $v_1$--torsion element of
$\ext^{2,j}(BP_*)$ with $j\equiv 0 \pmod 2$, 
then in the AAHSS for $\ext(BP_*P_{2m})$ the element
$z[2m]$ cannot be the target of a differential.
\end{prop}

\begin{proof}
The only elements which can support AAHSS differentials that hit $z[2m]$
are those of the form
$$ 1[2k-1] \qquad \text{or} \qquad \alpha_{i/j}[2k]. $$
We only need to consider elements in $t-s \equiv 1 \pmod 2$, since $z[2m]$
is in $t-2 \equiv 0 \pmod 2$.
The differentials given in Propositions~\ref{prop:AAHSSgendiffs} and 
\ref{prop:AAHSSv1diffs} tell us that these elements either kill or are
killed by other $v_1$--periodic elements.
\end{proof}

\begin{cor}\label{cor:betaroot}
The second filtered root invariant $\beta_{i/j} \in 
R^{[2]}_{BP}(\alpha_{i/j})$ always 
satisfies the second part of
\fullref{corA}.
\end{cor}

We now finish the proof of \fullref{thm:hroots} by verifying the second
part of \fullref{corA} each of our elements.

For the root invariants of the elements
$$ \wbar{\alpha}_{4/4},\quad  \wbar{\alpha}_{4/3},\quad \wbar{\alpha}_6 $$
we simply invoke \fullref{cor:betaroot}.  We mention that $\wbar{q}$ is
detected in the ANSS by the element $\beta_6$.

For the root invariants of the elements
$$ \wbar{\alpha}_{1}, \quad \wbar{\alpha}_1^2, \quad \wbar{\alpha}_1^3,
\quad
\wbar{\alpha}_{2/2}, \quad \wbar{\alpha}_{2}, \quad
\wbar{\alpha}_{4/4}\wbar{\alpha}_1, \quad
\wbar{\alpha}_{4/4}\wbar{\alpha}_1^2, \quad 
\wbar{\alpha}_{4/2}, \quad \wbar{\alpha}_{4}, \quad \wbar{\alpha}_{5},
\quad
\wbar{\alpha}_{6/2} $$
we look at the tables in \cite{Mahowaldmetastable} to see that the required
elements are non-zero in the homotopy of the stunted projective spaces.

The root invariant of $\wbar{\alpha}_{5}\wbar{\alpha}_1$ requires special
treatment.  We recall that we have the algebraic root invariant
$d_0^2 \in R_{\alg}(\wbar{\alpha}_{5}\wbar{\alpha}_1)$ carried by the $-19$--cell.  
We must verify that the
image of $\epsilon\wbar{\kappa} = \kappa^2$ under the map
$$ \pi_*(S^{-19}) \rightarrow \pi_*(P_{-19}) $$
is non-zero.  The problem is that when we look at Mahowald's computations
\cite{Mahowaldmetastable} we see that $d_0^2[-19]$ is killed in the AAHSS.
Indeed, we have the algebraic Atiyah--Hirzebruch differential
$$ d_6^{\mathrm{AAHSS}}(i[-13]) = d_0^2[-19]. $$
However, we have the Adams differential
$$ d_2^{\mathrm{ASS}}(i) = Pd_0h_0. $$
Therefore, in the $E_3$--term of the ASS for $P_{-19}$, 
the elements $Pd_0h_0[-13]$ and $Pd_0h_1[-14]$ have been equated, so we may
conclude that the combination of the AAHSS and ASS differentials implies
that the image of $\epsilon\wbar{\kappa}$ under the
inclusion of the bottom cell is detected by $Pd_0h_1[-14]$.  Mahowald's
tables \cite{Mahowaldmetastable} 
indicate that this element is non-zero in $\pi_*(P_{-19})$.

For the purposes of determining the root invariant of 
$$ \wbar{\alpha}_5\wbar{\alpha}_1^2 $$
we must determine whether the 
image of $\eta\wbar{q}$ under the map
$$ \pi_*(S^{-23}) \rightarrow \pi_*(P_{-23}) $$
is non-zero.  Examining the tables of Bruner \cite{Brunerweb}, we see that
the element $h_1q[-23]$ is non-trivial in $\ext(P_{-23})$.  We therefore
just need to check that it cannot be the target of a differential.  There
are three possible sources of a differential that would kill
$h_1q[-23]$.  These are represented by the elements 
$$ n[-20], \quad d_1[-21], \quad \text{and} \quad h_5h_2[-23].$$  
But these elements are permanent cycles, as argued in the following lemmas.

\begin{lem}
The element $n[-20] \in \ext^{5,16}(H_*P_{-23})$ is a permanent cycle.
\end{lem}

\begin{proof}
We just need to show that the element $\wbar{n} \in
\pi_{11}\bigl(S^{-20}\bigr)$
extends over $P_{-23}^{-20}$.  Since $2\wbar{n} = 0$ and $\eta \wbar{n} = 0$, 
it suffices to show that 
the Toda bracket
$$ \bra{2, \eta , \wbar{n}} $$
contains $0$.  
The element $\wbar{n}$ is given by
the Toda bracket
$$\wbar{n} \in \bra{\nu, \sigma, \wbar{\kappa}}$$
(see Mahowald--Tangora \cite{MahowaldTangora}).
We have
$$\bra{2, \eta, \wbar{n}}  = \bra{2, \eta, \bra{\nu, \sigma, \wbar{\kappa}}}
 \supseteq \bra{2, \eta, \nu, \sigma} \wbar{\kappa}.
$$
However, the Toda bracket $\bra{2,\eta, \nu, \sigma}$ lies in
$\pi_{13}(S^0)$, hence it must be zero.
\end{proof}

\begin{lem}
The element $d_1[-21] \in \ext^{4,15}(H_*P_{-23})$ is a permanent cycle.
\end{lem}

\begin{proof}
The element $\wbar{d}_1 \in \pi_{32}(S^0)$ extends over
$P_{-23}^{-21}$ to give an element which is detected by $d_1[-21]$, so we
may conclude that $d_1[-21]$ is a permanent cycle.
\end{proof}

The author thanks
W\,H~Lin for supplying the proof of the following lemma.

\begin{lem}
The element $h_5h_2[-23] \in \ext^{2,13}(H_*P_{-23})$ is a permanent cycle.
\end{lem}

\begin{proof}
The element $h_5h_2$ supports a differential of the form $d_3(h_5h_2) =
h_0p$ in the ASS for $S^0$.  However, in the AAHSS for
$\ext\bigl(H_*P_{-23}^{-22}\bigr)$, the
element $h_0p[-23]$ is killed by $p[-22]$.  The element $\wbar{p} \in
\pi_{33}(S^0)$ therefore extends to an element $\wbar{p}[-22] \in
\pi_{11}(P_{-23})$ which is detected by $h_5h_2[-23]$ in the ASS.  We
conclude that $h_5h_2[-23]$ is also a permanent cycle.  
\end{proof}

\bibliographystyle{gtart}
\bibliography{link}

\end{document}